\documentclass[twocolumn]{svjour3}  
\smartqed  % flush right qed marks, e.g., at end of proof

\usepackage[dvipdfmx]{color}
\usepackage{graphicx}
\usepackage[T1]{fontenc}
\usepackage{mathptmx}
\usepackage[scaled]{helvet}
\usepackage{courier}
\usepackage[subrefformat=parens]{subcaption}

\makeatletter
\let\cl@chapter\undefined
\makeatletter
\usepackage{amsmath,amssymb}   % for math commands
\usepackage{mathtools} % for math commands
\usepackage{cleveref}  % for clever reference
\Crefname{equation}{Eq.}{Eqs.}% 
\Crefname{figure}{Fig.}{Figs.}%
\usepackage{amsbsy}
\usepackage{bm}
\usepackage{algorithmic}
\usepackage{algorithm}
\usepackage{siunitx}
\usepackage{slashbox}

\begin{document}

\title{Review of the analytical prediction method of \\ surf-riding threshold in following sea,\\ and its relation to IMO second-generation intact stability criteria
%
%Practical methodology of wind process generation for examining autonomous operating systems of ships
}

\author{Atsuo Maki         \and Masahiro Sakai \and
        Tetsushi Ueta
}

\institute{Atsuo Maki \and Masahiro Sakaia \at
              Osaka University, 2-1 Yamadaoka, Suita, Osaka, Japan \\
              \email{maki@naoe.eng.osaka-u.ac.jp} %  \\
           \and
           Tetsushi Ueta \at
           Center for Administration of Information Technology, Tokushima University, 2-1 Minami-Josanjima-Cho, Tokushima 770-8506, Japan  \\
}

\date{Received: date / Accepted: date}

\maketitle

\begin{abstract}
In high-speed maritime operations, the broaching phenomenon can pose a significant risk when navigating in following/quartering seas. The occurrence of this phenomenon can result in a violent yaw motion, regardless of the steering effort, which, in turn, cause the resulting centrifugal force to capsize a vessel. A necessary condition for the occurrence of broaching is the surf-riding phenomenon. Therefore, the International Maritime Organization (IMO) has set up criteria to include theoretical formulas for estimating the occurrence of surf-riding phenomena. The theoretical equation used in the IMO's second-generation intact stability criteria (SGISC) to estimate the surf-riding threshold is based on Melnikov's method. This paper presents nonlinear equations describing the forward and backward motions of a ship. However, such equations cannot be directly solved; therefore, we proposed the use of and explain various approximate solution methods, including Meknikov's method. Subsequently, the relationship between the theoretical prediction method of the surf-riding threshold rooted in Melnikov's method and the IMO's SGISC is determined. 

\end{abstract}

\section{Introduction}
\begin{figure}[h]
        \centering
        \includegraphics[clip,width=1.0\hsize]{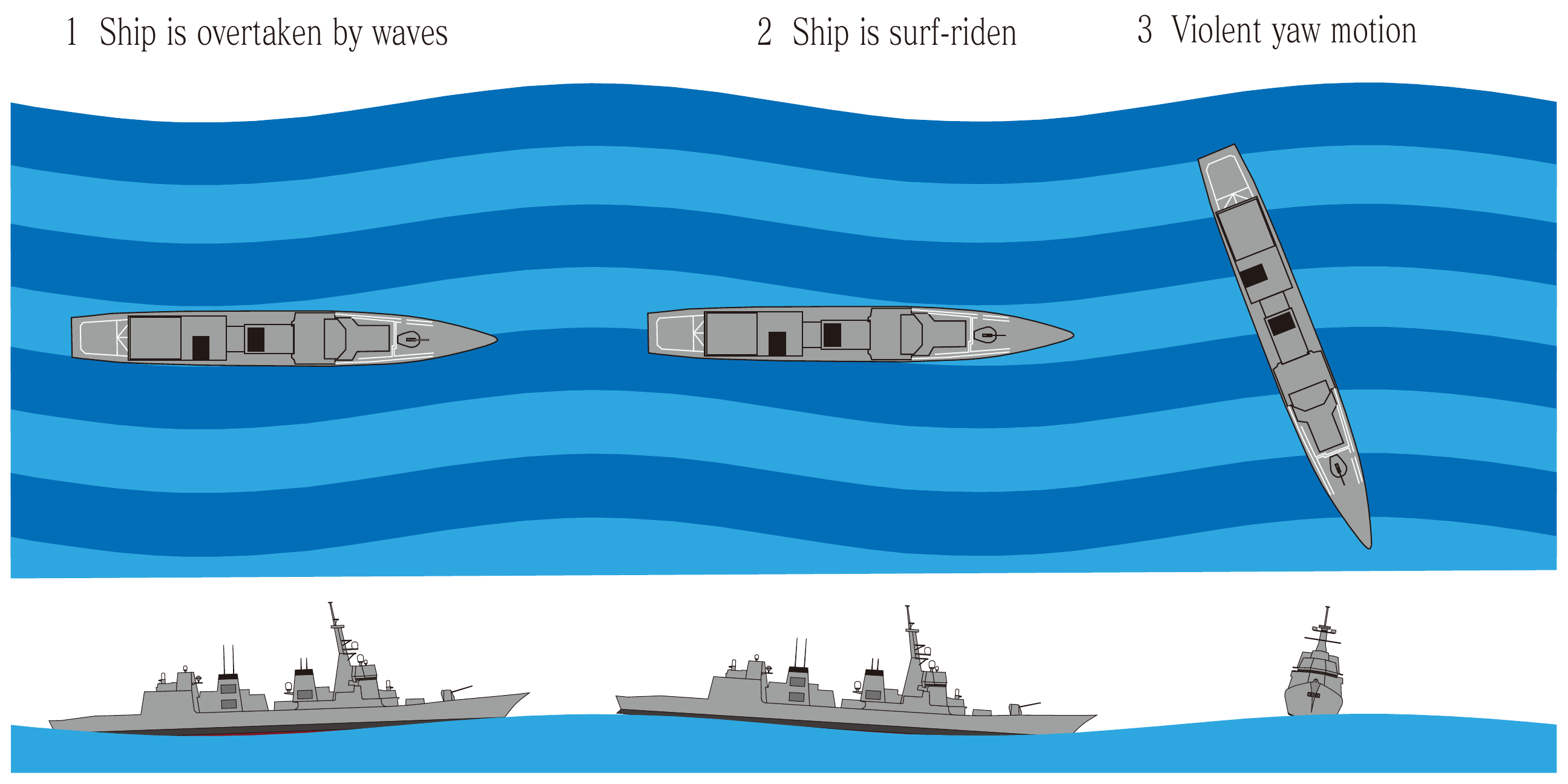}
        \caption{Schematic of the surf-riding phenomenon}
        \label{fig:Broaching_image}
    \end{figure}

    Marine vessels are known to operate at high speeds in following/quartering seas, where the broaching phenomenon may occur along with the surf-riding phenomenon. The broaching phenomenon is defined as a course-instability phenomenon among waves. As described in Fig.\ref{fig:Broaching_image}, initially, the ship is overtaken by waves; however, the force of the waves accelerates the hull, resulting in wave riding. On the wave down-slope, the course stability of the ship becomes unstable, resulting in a violent yaw motion. 
    %Broaching also has a path instability phenomenon that occurs at low speeds, which can be interpreted as a kind of self-excited oscillation, as shown by Spyrou~\cite{spyrou1996dynamic,spyrou1997dynamic}. Therefore, with appropriate control, the phenomenon of broaching at low speeds can be prevented. In the following, therefore, we will discuss the broaching phenomenon with wave riding that occurs when the ship is navigated at high speed.
    %
    Saunders described the broaching phenomenon that occurred with respect to the Portuguese Navy destroyer, N.R.P.~\textit{Lima}, in the Atlantic Ocean~\cite{saunders1965}. At the time, N.R.P.~\textit{Lima} was being operated at 26 knots in following waves with the wavelength-to-ship length ratio of almost 1. Despite the maximum steering efforts to maintain her course, she was forced to turn to port, and her roll angle was reported to have reached 67 degrees.

    Nicholson~\cite{nicholson1974} conducted an experiment on the model scale in the rectangular tank of the Hustler. A 5-m-long, twin-screw, twin-rudder model ship was navigated by manual steering with radio control in waves and a turning was recorded to the right even though the maximum rudder angle was 35 degrees to the left.

    Broaching describes the loss of ship maneuverability in waves and could occur even at low speeds, as described experimentally and theoretically by Kan~\cite{kan1990capsizing}, Spyrou~\cite{spyrou1996dynamic,spyrou1997dynamic}, and Maki~\cite{maki2009bifurcation}. However, this phenomenon is considered as a self-excited oscillation of the control system involving waves. Therefore, the accurate control can prevent the occurrence of broaching at low speeds. In this study, we focused on broaching at high speeds, which occurs after the occurrence of the surf-riding phenomenon. Such a type of broaching phenomenon has been extensively investigated from various aspects (e.g.~\cite{motora1982mechanism}). If a surf-riding phenomenon occurs, vessels cannot easily maneuver through even if the propeller thrust is reduced. In addition, a reduction of the propeller thrust results in a reduction in the rudder force. Therefore, from the standpoint of safety, the surf-riding phenomenon must be prevented.

    As shown later in text, the equations of ship motion in the following waves are the same as the equations of motion of a physical pendulum with respect to friction damping and constant motor torque. In the motion of such a pendulum, in the range of small torque, the motion is an oscillation that converges to an equilibrium point that exists, for example, at the lower point. The surfing phenomenon is identical to this situation. However, as the torque increases, the motion changes to a rotational motion, and this periodic state of motion is overtaken by waves. In terms of nonlinear dynamics, the boundary of such a motion is termed as the heteroclinic bifurcation. Several attempts have been made in the field of naval architect and ocean engineering to determine this bifurcation point.

    Grim~\cite{Grim1951} explained that under surf-riding threshold conditions, a trajectory leaving an unstable equilibrium point connects to another unstable equilibrium point. This phenomenon is known as a heteroclinic bifurcation. Makov~\cite{Makov1969} validated this finding by using the phase plane analysis. Annaniev~\cite{Ananiev1966} obtained the surf-riding threshold using the perturbation method. Kan~\cite{Kan_1990_surfriding} and Umeda~\cite{Umeda1990SR_irregular} also conducted the phase plane analyses, and based on these results, the International Maritime Organization (IMO) developed operational guidance for avoiding surf-riding in following seas as its MSC Circ.707~\cite{Imo1995}. Here, for simplicity, the critical Froude number for the surf-riding phenomenon was set as 0.3 for any ship. 
    
    The nonlinear surge equation, which has been used, cannot be solved analytically. Therefore, Kan applied Melnikov's method~\cite{holmes1980averaging} to propose an approximate formula for predicting the surf-riding threshold~\cite{Kan_1990_surfriding}. Spyrou~\cite{spyrou2006} extended Kan's approach and then obtained an approximate formula. Further, Maki generalized their approaches and then obtained an approximate formula~\cite{maki2010melnikov}. Spyrou~\cite{spyrou2001exact} also proposed an analytical formula for the nonlinear surge equation with a quadratically approximated damping component. Maki et al.~\cite{maki2010surfriding} applied the piecewise linear approximation to the sinusoidal wave force and obtained the analytical formula. Furthermore, Maki et al.~\cite{maki2014melnikov} applied the cubic polynomial approximation to the sinusoidal wave force to obtain the analytical formula. 
    
    The IMO has recently developed new-generation intact stability criteria, covering five failure modes, that is, parametric rolling, pure loss of stability, stability in dead-ship condition, broaching with surf-riding, and excessive accelerations. Accordingly, the risk of surf-riding phenomena must be assessed to prevent broaching phenomena. Therefore, Maki et al.~\cite{maki2010surfriding} proposed a calculation method directly based on Melnikov's method for this purpose. The current paper begins with a description of the nonlinear equations of motion that describe the surge motion in the waves. Then, several approximate solution methods of this nonlinear equation are described. They include the theoretical approximate formula described in IMO's new-generation intact stability criteria. Finally, the criteria are explained in detail.

\section{Notation}
    In the following equation, $\mathbb{R}$ represents a set of real numbers, and $\mathbb{R}^n$ denotes the $n$-dimensional Euclidean space. $\|x(t)\|$ for $x \in \mathbb{R}^n$ represents the Euclidean norm, $(x^T x)^{1/2}$. The absolute value of $x \in \mathbb{R}$ is denoted as $|x|$, where $i$ indicates the imaginary unit.
    
\section{Froude--Krylov surge force}
    As pointed out by Umeda~\cite{Umeda1983}, assuming that the hull form is almost longitudinally symmetric, the Froude-Krylov force is represented as a first-order approximation. Here, we show the detailed derivation of the sinusoidal surge force. 

    The velocity potential of incident wave, $\phi_0$, can be represented as
    \begin{equation}
        \begin{aligned}
            &\phi_0=-\frac{g \zeta_{\mathrm{W}}}{\omega} e^{-k_{\mathrm{W}} z_{\mathrm{s}}} \cdot\\
            &\sin \left(k_{\mathrm{W}} \xi_{\mathrm{G}}+k_{\mathrm{W}} x_{\mathrm{s}} \cos \chi_{\mathrm{s}} - k_{\mathrm{W}} y_{\mathrm{s}} \sin \chi_{\mathrm{s}} -\omega_{\mathrm{e}} t\right)
        \end{aligned}
    \end{equation}
    where $g$ is the gravitational acceleration, $\rho$ is the water density, $\zeta_{\mathrm{W}}$ is the amplitude of the incident wave, $k_{\mathrm{W}}=g/\omega^2=2 \pi / \lambda$ is the wave number, $\omega$ represents the wave frequency, $\lambda$ is the wavelength of the incident wave, and $\chi_{\mathrm{s}}$ is the ship direction to the incident wave.  
    Here, the following relation exists between wave and encounter frequencies:
    \begin{equation}
        \omega_{\mathrm{e}} = \omega - k_{\mathrm{W}} U \cos \chi_{\mathrm{s}}
    \end{equation}
    where $U$ is the forward velocity of a vessel. 
    The pressure can be calculated as
    \begin{equation}
        p=\rho\left(\frac{\partial}{\partial t}-U \frac{\partial}{\partial x_{\mathrm{s}}}\right) \phi_0
    \end{equation}
    In the following/quartering seas, the encounter frequency is almost zero. Therefore, hereafter, $\omega_{\mathrm{e}}=0$ and then $p$ can be calculated as
    \begin{equation}
        \begin{aligned}
            &p=-\rho g \zeta_{\mathrm{W}} e^{-k_{\mathrm{W}} z_{\mathrm{s}}} \cdot\\
            &\cos \left(k_{\mathrm{W}} \xi_{\mathrm{G}}+k_{\mathrm{W}} x_{\mathrm{s}} \cos \chi_{\mathrm{s}} -k_{\mathrm{W}} y_{\mathrm{s}} \sin \chi_{\mathrm{s}} \right)
        \end{aligned}
        \label{eq:pressue}
    \end{equation}
    The integration of pressure $p$ over the hull surface yields the Froude--Krylov force:
    \begin{equation}
        F_1^{\mathrm{FK}}=\iint_{S_{\mathrm{H}}}(-p) n_1 \mathrm{d} S
    \end{equation}
    The minus sign in $-p$, is derived from the definition of the normal vector.
    As the main goal here is to obtain the surge force, we need to have $x_{\mathrm{s}}$ directional normal vector $n_1$. However, unlike $y_{\mathrm{s}}$ or $z_{\mathrm{s}}$ directional normal vectors, the computation of $x_{\mathrm{s}}$ is not always easy. 
    Therefore, to bypass the use of $n_1$, Gauss's theorem is applied.

\begin{figure}[h]
        \centering
        \includegraphics[clip,width=0.9\hsize]{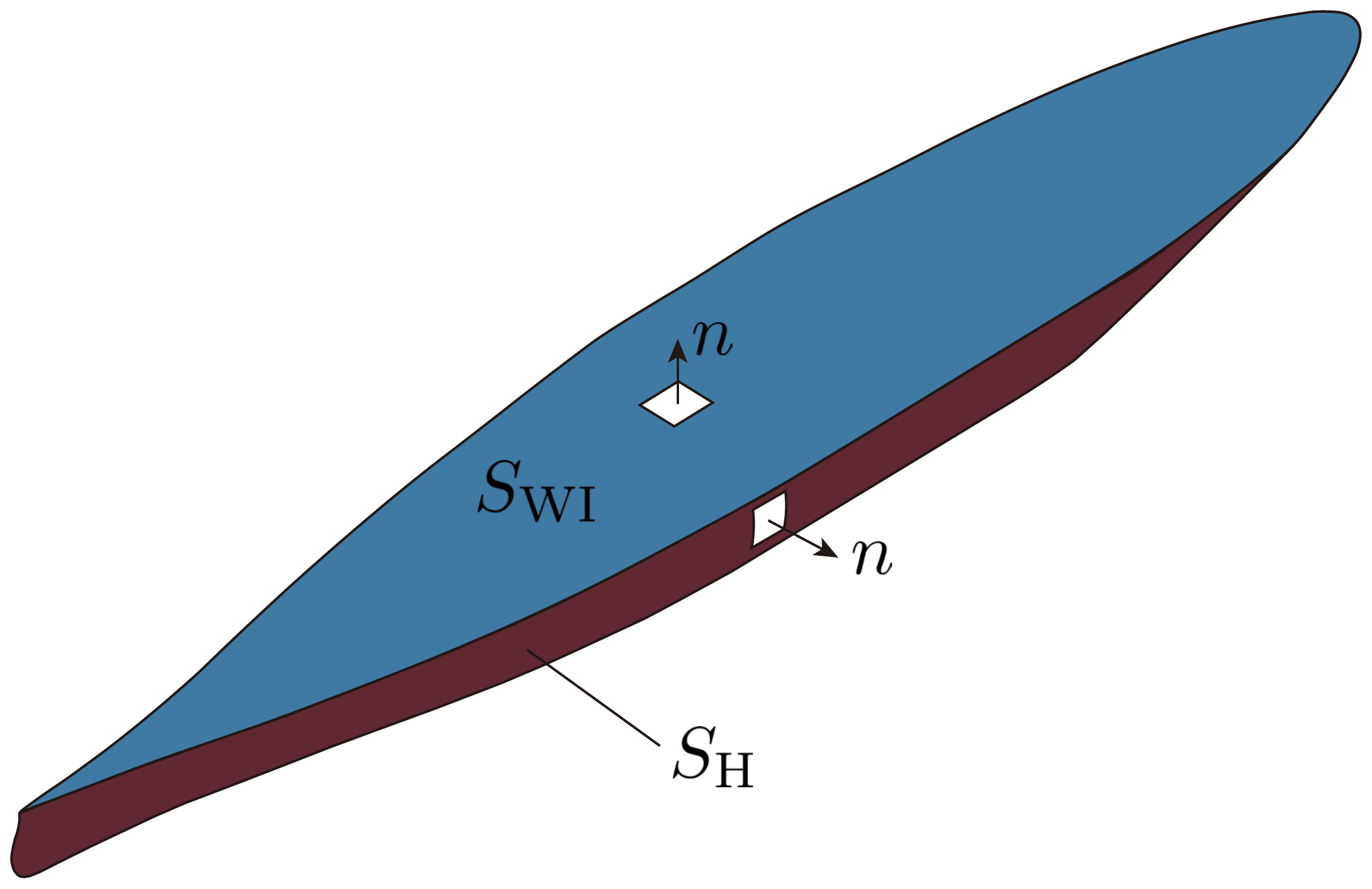}
        \caption{Definition of $S_{\mathrm{H}}$ and $S_{\mathrm{WI}}$}
        \label{fig:integral_domain}
    \end{figure}

    Here, as shown in Fig.~\ref{fig:integral_domain}, the authors introduce the surface region, $S_{\mathrm{WI}}$, i.e., the water-level region inside the hull, and the pressure on $S_{\mathrm{WI}}$ is rigorously zero in the framework of linear theory. 
    Then, by using Gauss's theorem, the surface integration is successfully converted to the volume integral. 
    \begin{equation}
        F_1^{\mathrm{FK}}=\iint_{S_{\mathrm{H}}+S_{\mathrm{WI}}}(-p) n_1 \mathrm{d} S=-\iiint_{V_0} \frac{\partial p}{\partial x_{\mathrm{s}}} \mathrm{d} V
        \label{eq:F1_Gauss}
    \end{equation}
    Now, we substitute Eq.~\ref{eq:pressue} into Eq.~\ref{eq:F1_Gauss}, and the volume integral can be divided into two integrals, namely $\mathrm{d}x$ and $\mathrm{d}s$, by following the scheme of the strip method.
    \begin{equation}
        \begin{aligned}
            &F_1^{\mathrm{FK}}  =-\rho g \zeta_{\mathrm{W}} k_{\mathrm{W}} \cos \chi_{\mathrm{s}} \cdot\\
            & \iiint_{V_0} e^{-k_{\mathrm{W}} z_{\mathrm{s}}} \sin \left(k_{\mathrm{W}} \xi_{\mathrm{G}} + k_{\mathrm{W}} x_{\mathrm{s}} \cos \chi_{\mathrm{s}} - k_{\mathrm{W}} y_{\mathrm{s}} \sin \chi_{\mathrm{s}} \right) \mathrm{d} V \\
            & =-\rho g \zeta_{\mathrm{W}} k_{\mathrm{W}} \cos \chi_{\mathrm{s}} \int_\mathcal{L} \mathrm{d} x \cdot\\
            & \iint_{S(x_{\mathrm{s}})} e^{-k_{\mathrm{W}} z_{\mathrm{s}}} \sin \left(k_{\mathrm{W}} \xi_{\mathrm{G}} + k_{\mathrm{W}} x_{\mathrm{s}} \cos \chi_{\mathrm{s}} - k_{\mathrm{W}} y_{\mathrm{s}} \sin \chi_{\mathrm{s}} \right) \mathrm{d} S
        \end{aligned}
    \end{equation}
    where $\mathcal{L}$ represents the x-directional integration over ship length $L$, and $S(x_{\mathrm{s}})$ represents the surface integral over a sectional area, $S(x_{\mathrm{s}})$.
    Now, the sinusoidal term can be expanded as
    \begin{equation}
        \begin{aligned}
            &\sin \left(k_{\mathrm{W}} \xi_{\mathrm{G}} + k_{\mathrm{W}} x_{\mathrm{s}} \cos \chi_{\mathrm{s}} - k_{\mathrm{W}} y_{\mathrm{s}} \sin \chi_{\mathrm{s}} \right) \\
            &=\sin \left(k_{\mathrm{W}} \xi_{\mathrm{G}} + k_{\mathrm{W}} x_{\mathrm{s}} \cos \chi_{\mathrm{s}} \right) \cos (k_{\mathrm{W}} y_{\mathrm{s}} \sin \chi_{\mathrm{s}})\\
            &-\cos \left(k_{\mathrm{W}} \xi_{\mathrm{G}} + k_{\mathrm{W}} x_{\mathrm{s}} \cos \chi_{\mathrm{s}} \right) \sin (k_{\mathrm{W}} y_{\mathrm{s}} \sin \chi_{\mathrm{s}})
        \end{aligned}
    \end{equation}
    Then, in the case of usual bilateral symmetry vessel, the second term in equation (8) disappears. Then, we can obtain
    \begin{equation}
        \begin{aligned}
            &F_1^{\mathrm{FK}}=-\rho g \zeta_{\mathrm{W}} k_{\mathrm{W}} \cos \chi_{\mathrm{s}} \cdot\\
            &\int_\mathcal{L} \sin \left(k_{\mathrm{W}} \xi_{\mathrm{G}} + k_{\mathrm{W}} x_{\mathrm{s}} \cos \chi_{\mathrm{s}} \right) \mathrm{d} x_{\mathrm{s}} \cdot\\
            &\iint_{S(x_{\mathrm{s}})} e^{-k_{\mathrm{W}} z_{\mathrm{s}}} \cos (k_{\mathrm{W}} y_{\mathrm{s}} \sin \chi_{\mathrm{s}}) \mathrm{d} S
        \end{aligned}
        \label{eq:F1_2nd_comp}
    \end{equation}
    Here, the authors assume the rectangular hull section, which has the same ship breadth, $B(x_{\mathrm{s}})$, and ship draft, $d(x_{\mathrm{s}})$. 
    Further, component $e^{-k_{\mathrm{W}} z_{\mathrm{s}}}$ is approximated as
    \begin{equation}
        e^{-k_{\mathrm{W}} z_{\mathrm{s}}} \approx e^{-k_{\mathrm{W}} d(x_{\mathrm{s}}) / 2}
    \end{equation}
    Then, the second integral in Eq.~\ref{eq:F1_2nd_comp} can be approximated as
    \begin{equation}
        \begin{aligned}
            &\iint_{S(x_{\mathrm{s}})} e^{-k_{\mathrm{W}} z_{\mathrm{s}}} \cos (k y_{\mathrm{s}} \sin \chi_{\mathrm{s}}) \mathrm{d} S \\
            &\approx e^{- k_{\mathrm{W}} d(x_{\mathrm{s}}) / 2} \int_0^{d(x_{\mathrm{s}})} \mathrm{d} z \int_{-B(x_{\mathrm{s}}) / 2}^{B(x_{\mathrm{s}}) / 2} \cos (k_{\mathrm{W}} y_{\mathrm{s}} \sin \chi_{\mathrm{s}}) \mathrm{d} y_{\mathrm{s}} \\
            & =C_1(x_{\mathrm{s}}) B(x_{\mathrm{s}}) \int_0^{d(x_{\mathrm{s}})} \mathrm{d} z_{\mathrm{s}} \\
            & \approx C_1(x_{\mathrm{s}}) S(x_{\mathrm{s}})
        \end{aligned}
    \end{equation}
    In equation (11), $C_1(x_{\mathrm{s}})$ represents
    \begin{equation}
        C_1(x_{\mathrm{s}}) \equiv \frac{\sin (k_{\mathrm{W}} \sin \chi B(x_{\mathrm{s}}) / 2)}{k_{\mathrm{W}} \sin \chi B(x_{\mathrm{s}}) / 2}
    \end{equation}
    Finally, the Froude--Krylov surge force can be written as
    \begin{equation}
        \begin{aligned}
            F_1^{\mathrm{FK}}=&-\rho g \zeta_{\mathrm{W}} k_{\mathrm{W}} \cos \chi \cdot \\
            &\int_\mathcal{L} e^{-k_{\mathrm{W}} d(x_{\mathrm{s}}) / 2} C_1(x_{\mathrm{s}}) S(x_{\mathrm{s}}) \cdot\\
            &\sin \left(k_{\mathrm{W}} \xi_{\mathrm{G}} + k_{\mathrm{W}} x_{\mathrm{s}} \cos \chi_{\mathrm{s}} \right) \mathrm{d} x_{\mathrm{s}}
        \end{aligned}
    \end{equation}
    In the case of following seas, $\chi_{\mathrm{s}}=0$; then, we have
    \begin{equation}
        \begin{aligned}
            &X_{\mathrm{W}}\left(\xi_{\mathrm{G}} \right)=-\rho g \zeta_{\mathrm{W}} k_{\mathrm{W}} \cdot\\
            &\int_{\mathrm{AE}}^{\mathrm{FE}} S(x_{\mathrm{s}}) e^{-k_{\mathrm{W}} d(x_{\mathrm{s}}) / 2} \sin k_{\mathrm{W}}\left(\xi_{\mathrm{G}}+x_{\mathrm{s}} \right) \mathrm{d} x_{\mathrm{s}}
        \end{aligned}
    \end{equation}
    Now, amplitude component $f$ and phase component $\varepsilon$ can be represented as
    \begin{equation}
        \begin{gathered}
            X_{\mathrm{W}}\left(\xi_{\mathrm{G}}\right)=-f_{\mathrm{W}} \sin \left(k_{\mathrm{W}} \xi_{\mathrm{G}} + \varepsilon \right)
        \end{gathered}
    \end{equation}
    where
    \begin{equation*}
        \left\{\begin{aligned}
                f_{\mathrm{W}} &= \sqrt{I_1^2+I_2^2} \\
                \varepsilon &= \tan ^{-1}\left(I_2 / I_1\right)\\
                I_1&=\rho g \zeta_{\mathrm{W}} k_{\mathrm{W}} \int_{A E}^{F E} S(x_{\mathrm{s}}) e^{-k_{\mathrm{W}} d(x_{\mathrm{s}}) / 2} \cos k_{\mathrm{W}} x_{\mathrm{s}} \mathrm{d} x_{\mathrm{s}} \\
                I_2&=\rho g \zeta_{\mathrm{W}} k_{\mathrm{W}} \int_{A E}^{F E} S(x_{\mathrm{s}}) e^{-k_{\mathrm{W}} d(x_{\mathrm{s}}) / 2} \sin k_{\mathrm{W}} x_{\mathrm{s}} \mathrm{d} x_{\mathrm{s}}
            \end{aligned}\right.
    \end{equation*}
    \color{black}
\section{Equation of motion}
    
\begin{figure*}
        \centering
        \includegraphics[clip,width=0.7\hsize]{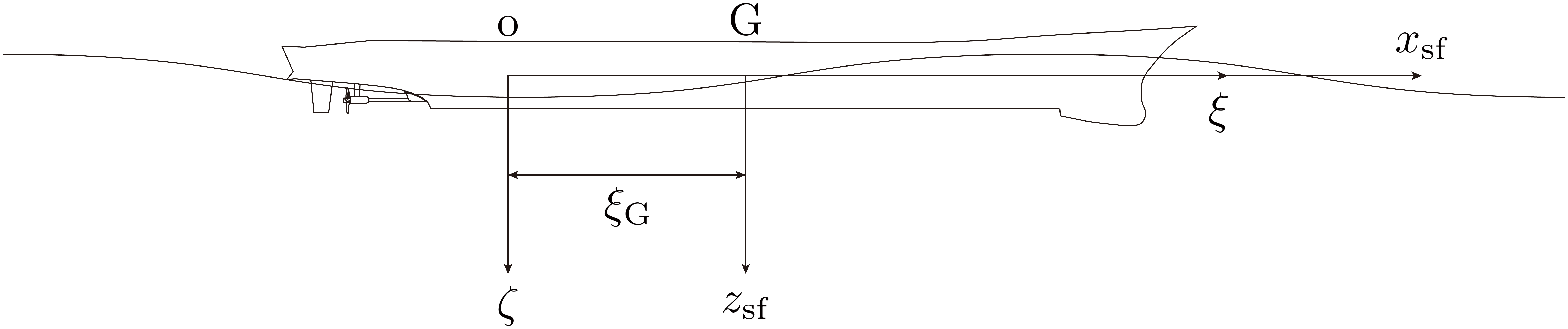}
        \caption{Coordinate systems}
        \label{fig:Coordinate_systems}
    \end{figure*}

    Two coordinate systems used in the paper are illustrated in Fig.~\ref{fig:Coordinate_systems}. The ship-fixed coordinate system, that is, $\mathrm{G}-x_{\mathrm{sf}}z_{\mathrm{sf}}$, has its origin at the center of gravity of the ship with the $x$ axis pointing toward the bow direction and $z$ axis pointing downward. An inertial coordinate system, namely $\mathrm{o}-\xi \zeta$, with the origin at a wave trough is employed with the $\xi$ axis pointing in the direction of wave travel and $\zeta$ axis pointing downward. 
    Therefore, the following relation on the velocity holds:
    \begin{equation}
        \frac{\mathrm{d} \xi_{\mathrm{G}}}{\mathrm{d} t} = c_{\mathrm{W}} + u
        \label{eq:relation_u_cw}
    \end{equation}
    Here, $c_{\mathrm{W}}$ represents the wave celerity. 
    
    The equation of motion that represents nonlinear surge motion is given as
    \begin{equation}
        \left(m+m_x\right) \frac{\mathrm{d}^2 \xi_{\mathrm{G}}}{\mathrm{d} t^2} + \left(R(u)-T(u, n_{\mathrm{P}})\right)-X_{\mathrm{W}}(\xi_{\mathrm{G}})=0
        \label{eq:equation_of_motion}
    \end{equation}
    In this equation, $R(u)$: the ship resistance which is positive in the negative x direction, $T(u, n_{\mathrm{P}})$: the propeller thrust, $m$: the ship mass, $m_x$: the added mass in the x direction.
    Hereafter, the resistance in waves is considered the same as the value in calm water at the same velocity; the effects of changes in trim, wetted surface, waterplane, local draughts, LCB position, and other factors are all ignored~\cite{Umeda1983}. 
    
    The ship resistance component, $R(u)$, is represented as n-th polynomial:
    \begin{equation}
        R(u) \approx \sum_{j=1}^{n_{\mathrm{M}}} r_j u^j = r_1 u + r_2 u^2 + \cdots,
    \end{equation}
    Then, the ship thrust can be calculated as
    \begin{equation}
        \begin{gathered}
            T_{\mathrm{e}}(u, n_{\mathrm{P}}) = \left(1-t_{\mathrm{P}}\right) \rho n^2 D_{\mathrm{P}}^4 K_{\mathrm{T}}\left(J( u, n_{\mathrm{P}} ) \right)
        \end{gathered}
    \end{equation}
    where
    \begin{equation*}
        J ( u, n_{\mathrm{P}} )=\frac{(1-w_{\mathrm{P}}) u}{n_{\mathrm{P}} D_{\mathrm{P}}}
    \end{equation*}
    where $t_{\mathrm{P}}$ and $w_{\mathrm{P}}$ are the thrust deduction and wake fraction. Thrust coefficient $K_{\mathrm{T}}\left(J(u, n_{\mathrm{P}})\right)$ inside $T_{\mathrm{e}}$ can be calculated as 
    \begin{equation}
        \begin{aligned}
            &K_{\mathrm{T}}(J( u, n_{\mathrm{P}} )) \approx \sum_{j=0}^{n_{\mathrm{M}}} \kappa_i J^j( u, n_{\mathrm{P}} ) \\
            &= \kappa_0+\kappa_1 J( u, n_{\mathrm{P}} )+\kappa_2 J^2( u, n_{\mathrm{P}} )+\cdots
        \end{aligned}
    \end{equation}
    By substituting Eq.\ref{eq:relation_u_cw} in Eq. (20), the ship thrust can be calculated as
    \begin{equation}
        \begin{aligned}
            T_{\mathrm{e}}(u, n_{\mathrm{P}}) = \sum_{j=0}^{n_{\mathrm{M}}} \frac{\left(1-t_{\mathrm{P}}\right)\left(1-w_{\mathrm{P}}\right)^j \rho \kappa_j u^j}{n_{\mathrm{P}}^{j-2} D_{\mathrm{P}}^{j-4}},
        \end{aligned}
    \end{equation}
    Then, the equation of motion becomes
    \begin{equation}
        \begin{aligned}
            &(m +m_x) \frac{\mathrm{d}^2 \xi_{\mathrm{G}}}{\mathrm{d} t^2} + \sum_{j=1}^{n_{\mathrm{M}}} \sum_{k=1}^j c_j(n_{\mathrm{P}}) \left(\begin{array}{l}
            j \\
            k
            \end{array}\right) \left( \frac{\mathrm{d} \xi_{\mathrm{G}}}{\mathrm{d} t} \right)^k c_{\mathrm{W}}^{j-k}\\
            &+f_{\mathrm{W}} \sin k_{\mathrm{W}} \xi_{\mathrm{G}} =T_{\mathrm{e}}\left(c_{\mathrm{W}}, n_{\mathrm{P}} \right)-R\left(c_{\mathrm{W}}\right)
        \end{aligned}
    \end{equation}
    where
    \begin{equation*}
        c_j(n_{\mathrm{P}}) \equiv r_j - \frac{\left(1-t_{\mathrm{P}}\right)\left(1-w_{\mathrm{P}}\right)^j \rho \kappa_j}{n_{\mathrm{P}}^{j-2} D^{j-4}}
    \end{equation*}
    The introduction of new coefficients yields the following results:
    \begin{equation}
        \begin{aligned}
            &\frac{\mathrm{d}^2 \xi_{\mathrm{G}}}{\mathrm{d} t^2} + \sum_{k=1}^{n_{\mathrm{M}}} a_k(n_{\mathrm{P}}) \left( \frac{\mathrm{d} \xi_{\mathrm{G}}}{\mathrm{d} t} \right)^k \\
            &+ \frac{f_{\mathrm{W}}}{m+m_x} \sin k_{\mathrm{W}} \xi_{\mathrm{G}} = \frac{T_{\mathrm{e}}\left(c_{\mathrm{W}}, n_{\mathrm{P}} \right)-R\left(c_{\mathrm{W}}\right)}{m + m_x}
        \end{aligned}
        \label{eq:equation_of_motion_for_use}
    \end{equation}
    where
    \begin{equation*}
        \begin{gathered}
             a_k(n_{\mathrm{P}}) = \frac{1}{m + m_x} \sum_{j=k}^{n_{\mathrm{M}}} c_j(n_{\mathrm{P}}) \left(\begin{array}{l}
                j \\
                k
            \end{array}\right) 
            c_{\mathrm{W}}^{j-k}
        \end{gathered}
    \end{equation*}
    The final obtained equation of motion (Eq.~\ref{eq:equation_of_motion_for_use}) is identical to that of a pendulum swing with constant torque.

    Now, the introduction of new coefficients yields the following:
    \begin{equation}
        \begin{gathered}
            \frac{\mathrm{d}^2 y}{\mathrm{d} t^2} + \sum_{k=1}^{n_{\mathrm{M}}} A_k(n_{\mathrm{P}}) \left( \frac{\mathrm{d} y}{\mathrm{d} t} \right)^k + q \sin y = r\left(n_{\mathrm{P}}\right)
        \end{gathered}
        \label{eq:eq:equation_of_motion_nondim_1}
    \end{equation}
    where
    \begin{equation*}
        \left\{\begin{aligned}
            y &= k \xi_G\\
            A_k(n_{\mathrm{P}}) &= \frac{1}{k_{\mathrm{W}}^{k-1} (m + m_x)} \sum_{j=k}^{n_{\mathrm{M}}} c_j(n_{\mathrm{P}}) \left(\begin{array}{l}
                j \\
                k
            \end{array}\right) 
            c_{\mathrm{W}}^{j-k}\\
            q &= \frac{f_{\mathrm{W}}k_{\mathrm{W}}}{m + m_x}\\
            r\left(n_{\mathrm{P}}\right) &= \frac{\left(T_{\mathrm{e}}\left(n_{\mathrm{P}}\right)-R\right) k_{\mathrm{W}}}{m+m_x}
        \end{aligned}\right.
    \end{equation*}

    Further, if the following nondimensional time $\tau = \sqrt{q} t$ is introduced, then coefficient $q$ in front of sinusoidal term can be eliminated:
    \begin{equation}
        \begin{gathered}
            \frac{\mathrm{d}^2 y}{\mathrm{d} \tau^2} + \sum_{k=1}^{n_{\mathrm{M}}} \bar{A}_k(n_{\mathrm{P}}) \left( \frac{\mathrm{d} y}{\mathrm{d} \tau} \right)^k + \sin y = \bar{r}\left(n_{\mathrm{P}}\right)
        \end{gathered}
        \label{eq:equation_of_motion_nondim_2}
    \end{equation}
    where
    \begin{equation*}
        \left\{\begin{aligned}
            \tau &= \sqrt{q}t\\
            \bar{A}_k(n_{\mathrm{P}}) &= A_k(n_{\mathrm{P}}) q^{k/2-1}\\ 
            \bar{r}\left(n_{\mathrm{P}}\right) &= \frac{r\left(n_{\mathrm{P}}\right)}{ q}
        \end{aligned}\right.
    \end{equation*}
    This equation has two nonlinear components, complicating the direct analytical approach: (1) the ``damping'' component, $\frac{\mathrm{d}y}{\mathrm{d} \tau}$ and (2) the ``restoring'' component, $y$.

\begin{figure}[htb]
        \centering
        \includegraphics[clip,width=1.0\hsize]{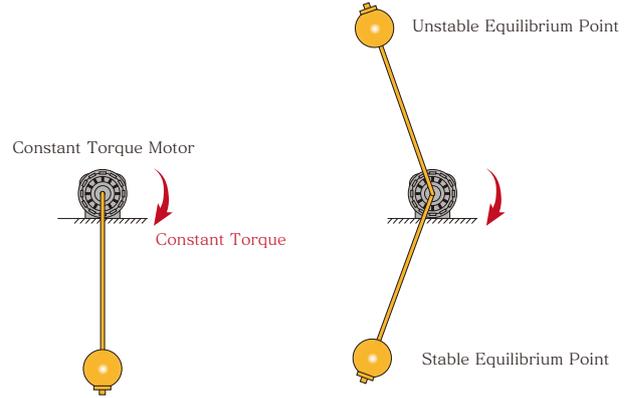}
        \caption{Equivalent pendulum with the surf-riding motion.}
        \label{fig:pendulum}
    \end{figure}

\begin{figure}[htb]
        \centering
        \includegraphics[clip,width=1.0\hsize]{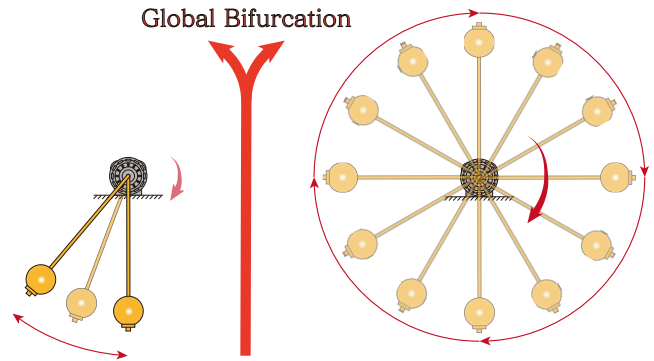}
        \caption{Global bifurcation in the pendulum }
        \label{fig:pendulum_bifurcation}
    \end{figure}

    The obtained equation is equivalent to the equation of motion of a nonlinear pendulum with a constant-torque motor, as shown in the left panel of Fig.~\ref{fig:pendulum}. This equation yields two equilibrium points on the lower and upper sides, as shown in the right panel of Fig.~\ref{fig:pendulum}, which are stable and unstable, respectively. As shown in Fig.~\ref{fig:pendulum_bifurcation}, when the constant torque is small, the pendulum oscillates around the lower, stable equilibrium point, as shown in the left panel. However, when the magnitude of the torque is increased, this pendulum changes to a rotational motion, moving in circles, as shown in the right panel. The boundary between the two is the global bifurcation point, which in this case is called the homoclinic bifurcation point.

\section{Phase portrait}
    \label{section_PhasePortrait}
\begin{figure}[htb]
        \centering
        \includegraphics[clip,width=1.0\hsize]{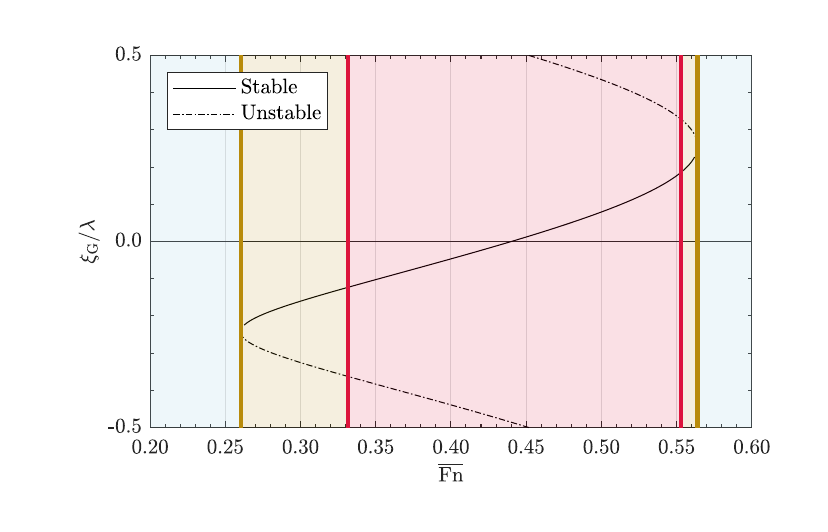}
        \caption{Position of the equilibrium point and its stability with respect to $\lambda/L=1.25$ and $H/\lambda=0.04$}
        \label{fig:Equilibrium_point}
    \end{figure}

    First, the longitudinal position of the equilibrium point is plotted in Fig.~\ref{fig:Equilibrium_point}, showing the ship's longitudinal position as $\xi_{\mathrm{G}}/\lambda$, i.e., $y$, which satisfies the following equation:
    \begin{equation}
        q \sin y - r\left( n_{\mathrm{P}} \right) = 0
    \end{equation}
    The horizontal axis of this plot displays the Froude number, $\mathrm{\overline{Fn}}$, which corresponds to the Froude number for the same propeller revolution in calm water. In addition, the vertical axis is a longitudinal position. Here, $\mathrm{\overline{Fn}}$ is nondimensional ship speed, and is defined as
    \begin{equation*}
        \mathrm{\overline{Fn}} = \frac{u}{\sqrt{Lg}}
    \end{equation*}
    As no equilibrium point is observed in the blue region, a ship cannot surf-ride. The boundary between blue and yellow regions is the local bifurcation point, which is known as the tangent bifurcation, and the tangent bifurcation points are located at $\mathrm{\overline{Fn}}=0.2602$ and $\mathrm{\overline{Fn}}=0.5639$. However, note that the surf-riding phenomenon does not always occur in yellow areas; it could occur if special initial values are used~\cite{umeda1990SR_regular,Umeda1990SR_irregular}. The boundary between yellow and red regions is the global bifurcation point, termed as the heteroclinic bifurcation, and the corresponding heteroclinic bifurcation points are located at $\mathrm{\overline{Fn}}=0.3318$ and $\mathrm{\overline{Fn}}=0.5532$. In addition, surf-riding occurs in this red region for all of the initial values. Hereafter, the heteroclinic bifurcation point at lower velocity is referred to as the \textit{\textbf{surf-riding threshold}} and the heteroclinic bifurcation point at high velocities is termed as the \textit{\textbf{wave-blocking threshold}}.

\begin{figure*}
        \centering
        \includegraphics[clip,width=\hsize]{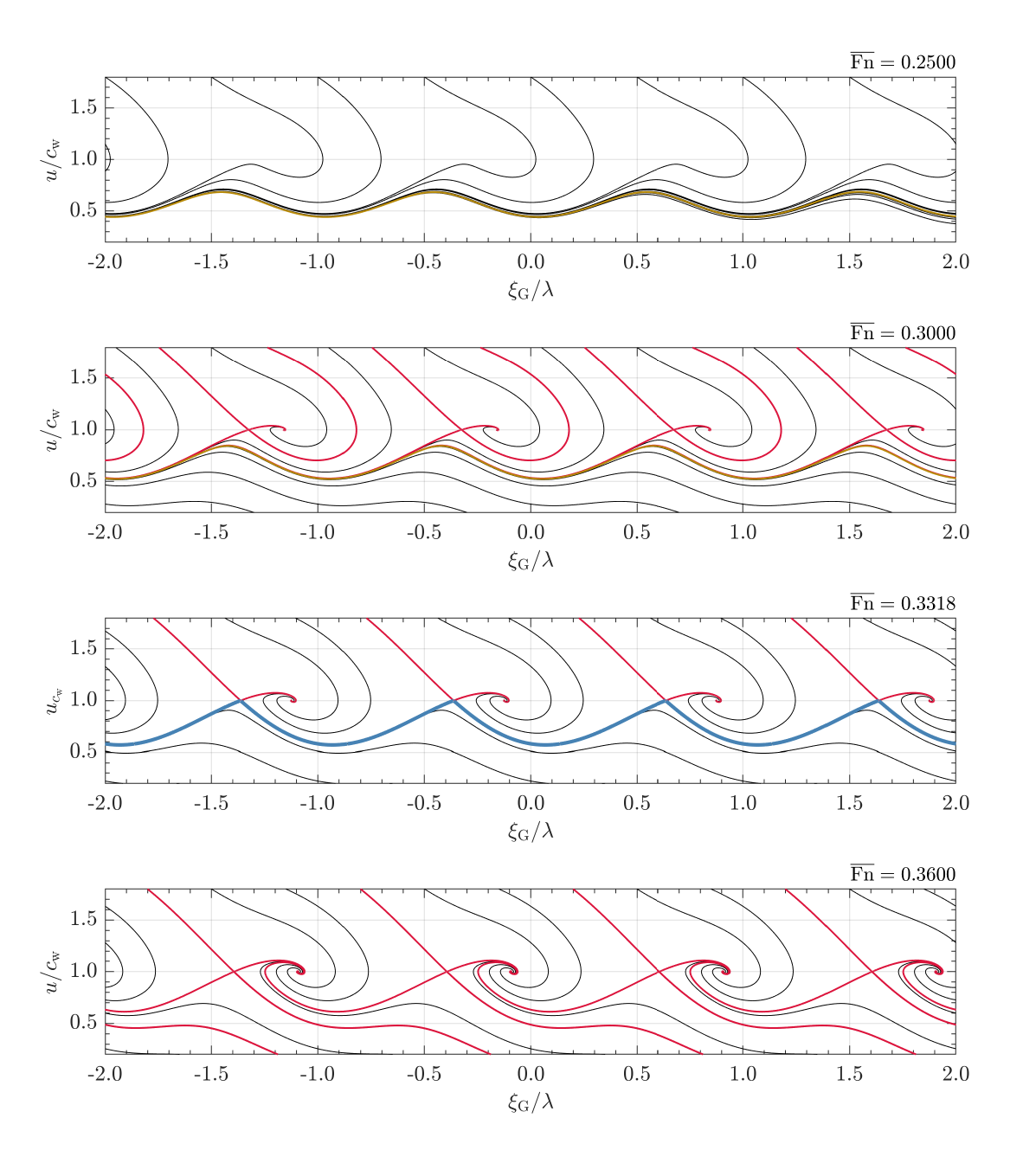}
        \caption{Phase portrait in the case of low speed with $\lambda/L=1.25$ and $H/\lambda=0.04$. Here, the gold lines represent the periodic attractors, red lines represent stable and unstable trajectories, and blue lines represent the heteroclinic trajectories.}
        \label{fig:Phase_portrait_lower}
    \end{figure*}

\begin{figure*}
        \centering
        \includegraphics[clip,width=\hsize]{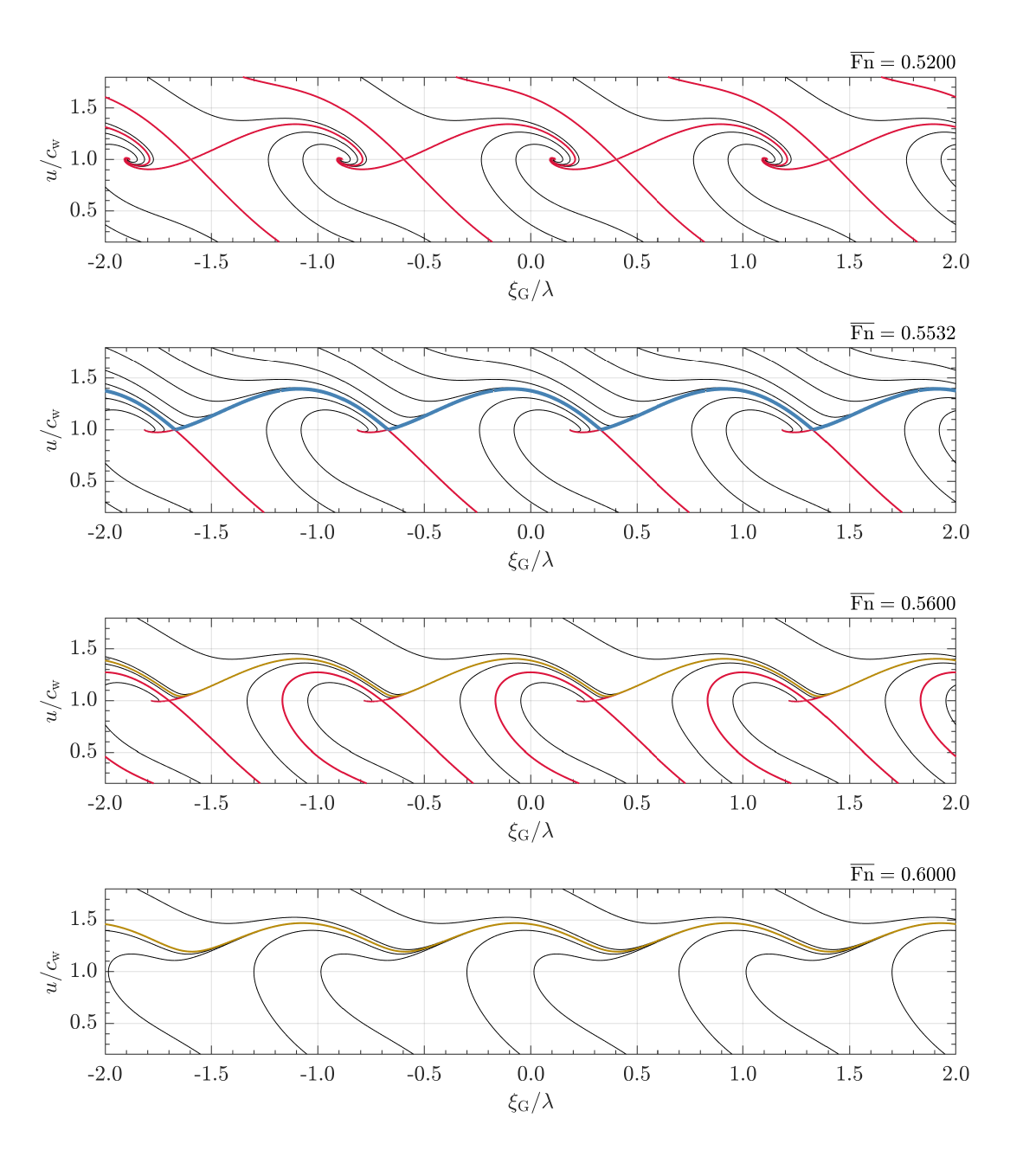}
        \caption{Phase portrait in the case of high speed with $\lambda/L=1.25$ and $H/\lambda=0.04$. Here, the gold lines represent the periodic attractors, red lines represent the stable and unstable trajectories, and blue lines represent the heteroclinic trajectories.}
        \label{fig:Phase_portrait_upper}
    \end{figure*}

    The phase portraits for DTMB5415 are shown in Figs.~\ref{fig:Phase_portrait_lower} and \ref{fig:Phase_portrait_upper}. Details of the subject vessel used are shown in Section~\ref{sec:comparison}. The wave conditions in this calculation are $\lambda/L=1.25$ and $H/\lambda=0.04$. Fig.~\ref{fig:Phase_portrait_lower} indicates the change in the phase portrait due to the change in the ship speed around the surf-riding threshold. With $\overline{\mathrm{Fn}}=0.2500$, every trajectory converges to the periodic attractor (gold line), as no stable equilibrium point exists for this condition. With $\overline{\mathrm{Fn}}=0.3300$, there exist trajectories that converge to the periodic attractor and stable equilibrium point. The stable equilibrium point indicates the surf-riding phenomenon, and this attractor has been analyzed in detail by Umeda~\cite{umeda1990SR_regular,Umeda1990SR_irregular}. In this wave condition, the surf-riding threshold is $\overline{\mathrm{Fn}}=0.3318$, and the blue line existing in the third plot indicates the heteroclinic orbit connecting the saddles. For $\overline{\mathrm{Fn}}=0.3600$, every trajectory converges to a stable equilibrium point. The phase plane changes qualitatively before and after the surf-riding threshold. Fig.\ref{fig:Phase_portrait_upper} shows the change in the phase portrait due to the change in the ship speed around the wave-blocking threshold. In this wave condition, the wave-blocking threshold is $\overline{\mathrm{Fn}}=0.5532$. Similar to the qualitative change in the phase plane near the surf-riding threshold, a qualitative change occurs in the phase plane around the wave-blocking threshold.
    

\section{Numerical bifurcation analysis~\cite{maki2009bifurcation}}
    Numerical techniques for directly identifying heteroclinic bifurcation points have been extensively investigated ~\cite{kuznetsov1998elements}. This section describes the numerical method proposed by Kawakami ~\cite{kawakami1981separatrix}. Although not employed by Maki~\cite{maki2009bifurcation}, the numerical bifurcation tool, AUTO, developed by Doedel and Friedman\cite{doedel1990numerical,friedman1991numerical} has been widely used in related studies. 
    
    First, the state vector is defined as
    \begin{equation}
        \mathcal{X} = \left( \xi_{\mathrm{G}} / \lambda, u \right)^{\top} \in \mathbb{R}^2
    \end{equation}
    The equation of motion, i.e., Eq.~\ref{eq:equation_of_motion_for_use}, is represented as the following vector form:
    \begin{equation}
            \frac{\mathrm{d} \mathcal{X}}{\mathrm{d} t} = f(\mathcal{X}) ~~\left(= (f_1(\mathcal{X}), f_2(\mathcal{X}))^{\top} \right)\
    \end{equation}
    where
    \begin{equation*}
        \left\{ \begin{aligned}
                f_1(\mathcal{X}) &= u / \lambda\\
                f_2(\mathcal{X}) &= - \frac{R(u)-T(u, n_{\mathrm{P}})}{m + m_x} + \frac{X_{\mathrm{W}}(\xi_{\mathrm{G}})}{m + m_x}=0
            \end{aligned}\right.
    \end{equation*}
    Here, $f:\mathbb{R}^2 \rightarrow \mathbb{R}^2$, $f_{1, 2}:\mathbb{R}^2 \rightarrow \mathbb{R}$.
    
    The definitions of unstable manifold $\alpha\left(x_0 \right)$ and stable manifold $\omega\left(x_1 \right)$ are given as follows:
    \begin{equation}
        \left\{\begin{aligned}
            \alpha\left(\mathcal{X}_0\right)&=\left\{\mathcal{X} \mid \lim _{t \rightarrow-\infty} \mathcal{X}(t) = \mathcal{X}_0\right\}  \quad \\
            \quad \omega\left(\mathcal{X}_1\right)&=\left\{\mathcal{X} \mid \lim _{t \rightarrow+\infty} \mathcal{X}(t) = \mathcal{X}_1\right\}
        \end{aligned}\right.
    \end{equation}
    The saddle-type equilibrium points are defined as $\mathcal{X}_0 \equiv\left(\xi_{\mathrm{G}0}/\lambda, c_{\mathrm{W}}\right)^{\top}$ and $\mathcal{X}_1 \equiv\left(\xi_{\mathrm{G}1}/\lambda, c_{\mathrm{W}}\right)^{\top}$, and $\xi_{\mathrm{G}0}/\lambda = \xi_{\mathrm{G}1}/\lambda \pm 2\pi$.
    Here, $\mathcal{X}_0$ and $\mathcal{X}_1$ satisfy
    \begin{equation}
        f\left(\mathcal{X}_0, n_{\mathrm{P}} \right) = f\left(\mathcal{X}_1 ; n_{\mathrm{P}} \right)=0
        \label{eq:Newton_cond_1}
    \end{equation}. Now, we linearize the state equation in the vicinity of $\mathcal{X}_0$ and $\mathcal{X}_1$ as follows:
    \begin{equation}
        \frac{\mathrm{d}\mathcal{X}}{\mathrm{d} t} = A \cdot\left(\mathcal{X}-\mathcal{X}_0\right) \quad \text{and} \quad \frac{\mathrm{d} \mathcal{X}}{\mathrm{d} t} =A \cdot \left(\mathcal{X}-\mathcal{X}_1\right)
    \end{equation}
    $A \in \mathbb{R}^{2 \times 2}$ indicates the Jacobi matrix at $x_0$ or $x_1$ as follows:
    \begin{equation}
        A = \left(\begin{array}{ll}
        \alpha_{11} & \alpha_{12} \\
        \alpha_{21} & \alpha_{22}
        \end{array}\right) = \left. \left(\begin{array}{cc}
        0 & 1 / \lambda \\
        \frac{\partial f_2}{\partial \xi_{\mathrm{G}/\lambda}} & \frac{\partial f_2}{\partial u}
        \end{array}\right) \right|_{\mathcal{X}=\mathcal{X}_0 \;\; \text{or} \;\; \mathcal{X}_1}
    \end{equation}
    Here, $\mu_\alpha$ and $\mu_\omega$ represent eigenvalues of $A$, which represent the negative and positive real parts, respectively:
    \begin{equation}
        \begin{aligned}
            &2 \mu_{\alpha,\omega}-\left(\alpha_{11}+\alpha_{22}\right)\\
            &\mp \sqrt{\left(\alpha_{11}+\alpha_{22}\right)^2-4\left(\alpha_{11} \alpha_{22}-\alpha_{12} \alpha_{21}\right)}=0
        \end{aligned}
        \label{eq:Newton_cond_2}
    \end{equation}
    The eigenvectors for these two eigenvalues are defined as $h_\alpha \equiv\left(\xi_{\mathrm{G}\alpha}/\lambda, u_\alpha\right)^{\top}$ and $h_\omega \equiv\left(\xi_{\mathrm{G}\omega}/\lambda, u_\omega\right)^{\top}$. Now, the authors take $\mathcal{X}_\alpha$ and $\mathcal{X}_\omega$ on $\alpha$ and $\omega$ blanches as follows:
    \begin{equation}
        h_\alpha = \mathcal{X}_\alpha - \mathcal{X}_0 \quad \text{and} \quad h_\omega = \mathcal{X}_\omega - \mathcal{X}_1
    \end{equation}
    Of course, $h_\alpha$ and $h_\omega$ satisfy the following condition of eigen direction:
    \begin{equation}
        \left(A - \mu_\alpha I\right) \cdot h_\alpha = 0 \quad \text{and} \quad \left(A - \mu_\omega I\right) \cdot h_\omega = 0
        \label{eq:Newton_cond_3}
    \end{equation}
    Here, $I \in \mathbb{R}^2$ indicates the unit matrix. Further, $\mathcal{X}_{\alpha}$ and $\mathcal{X}_{\omega}$ are located in the vicinity of the saddles; therefore, the authors impose the following condition:
    \begin{equation}
        \|h_{\alpha}\|^2 = \|h_{\omega}\|^2 = \varepsilon_{\mathrm{h}}^2
        \label{eq:Newton_cond_4}
    \end{equation}
    Here, $\varepsilon_{\mathrm{h}} \ll 1$ is a preset sufficiently small value. 

    Moreover, heteroclinic bifurcation occurs when stable and unstable manifolds contact each other. Now, we define trajectory $\psi\left(\mathcal{X}_{\mathrm{I}}, \tau_{\mathrm{I}}\right) \in \mathbb{R}^2$ at time $\tau_{\mathrm{I}}$ with initial condition $\mathcal{X}_{\mathrm{I}}$ as \begin{equation}
        \psi\left(\mathcal{X}_{\mathrm{I}}, \tau_{\mathrm{I}} \right) = (\psi_1\left(\mathcal{X}_{\mathrm{I}}, \tau_{\mathrm{I}}\right), \psi_2\left(\mathcal{X}_{\mathrm{I}}, \tau_{\mathrm{I}} \right))^{\top}
    \end{equation}
   Trajectories starting from $\mathcal{X}_0$ and $\mathcal{X}_1$ contact at the intermediate point at time $\tau_{\mathrm{I}}$.
   \begin{equation}
       \left\{\begin{array}{l}
        \varphi\left(\mathcal{X}_\alpha, \tau_{\mathrm{I}} \right)-\varphi\left(\mathcal{X}_\omega, -\tau_{\mathrm{I}} \right)=0 \\
        \psi\left(\mathcal{X}_\alpha, \tau_{\mathrm{I}} \right)-\psi\left(\mathcal{X}_\omega, -\tau_{\mathrm{I}} \right)=0
        \end{array}\right.
        \label{eq:Newton_cond_5}
   \end{equation}
    The computation of $\psi\left(\mathcal{X}_{\omega}, -\tau_{\mathrm{I}} \right)$ is obtained by solving ODE in inverse time. 
    
    Finally, every condition (Eq.~\ref{eq:Newton_cond_1}, Eq.~\ref{eq:Newton_cond_2}, Eq.~\ref{eq:Newton_cond_3}, Eq.~\ref{eq:Newton_cond_4}, Eq.~\ref{eq:Newton_cond_5}) was solved by using Newton's method for $n_{\mathrm{P}}$, $\mathcal{X}_0$, $\mathcal{X}_1$, $\mu_{\alpha}$, $\mu_{\omega}$, $h_{\alpha}$, $h_{\omega}$, and $\tau_{\mathrm{I}}$.

\section{Method based on the quadratic approximation of the damping component ~\cite{spyrou2006}}\label{sec:Spyrou}

    As stated previously, the equation of motion (Eq.~\ref{eq:equation_of_motion}) has two components that complicate the theoretical approach. Spyrou~\cite{spyrou2006} approximated the nonlinear ``damping'' component by using the quadratic function as follows:
    
    \begin{equation}
        \begin{gathered}
            \frac{\mathrm{d}^2 y}{\mathrm{d} \tau^2} + \gamma(n_{\mathrm{P}}) \operatorname{sgn} \left( \frac{\mathrm{d} y}{\mathrm{d} \tau} \right) \cdot \left( \frac{\mathrm{d} y}{\mathrm{d} \tau} \right)^2 + \sin y = \bar{r}\left(n_{\mathrm{P}}\right)
        \end{gathered}
        \label{eq:equation_of_motion_spyrou}
    \end{equation}
    where
    \begin{equation*}
        \gamma(n_{\mathrm{P}}) = -\frac{\displaystyle \sum_{j=1}^{n_{\mathrm{M}}} \bar{A}_j \displaystyle \int_0^{v_\mathrm{e}} v^{j + 2} \mathrm{d}v }{\displaystyle \int_0^{v_\mathrm{e}} v^4 \mathrm{d} v}
    \end{equation*}
    where $v = \mathrm{d} y / \mathrm{d} \tau$ and $v_\mathrm{e}$ is an upper limit of the least square fit of the damping component. In addition, Eq.~\ref{eq:equation_of_motion_spyrou} is well known to possess an analytical solution. 
    Furthermore, the variable transformation is applied to Eq.~\ref{eq:equation_of_motion_spyrou} to yield
    \begin{equation}
        \frac{1}{2} \frac{\mathrm{d} v^2}{\mathrm{d} y}+\gamma(n_{\mathrm{P}}) \operatorname{sgn} v \cdot v^2=-\sin y + \bar{r}\left(n_{\mathrm{P}}\right)
    \end{equation}
    In the case of $\frac{\mathrm{d}y}{\mathrm{d} \tau}<0$, the solution becomes
    \begin{equation}
        \frac{\mathrm{d}y}{\mathrm{d} \tau} = -\sqrt{c_2 e^{2 \gamma(n_{\mathrm{P}}) y}+\frac{2(\cos y+2 \gamma(n_{\mathrm{P}}) \sin y)}{1+4 \gamma^2(n_{\mathrm{P}})}-\frac{\bar{r}\left(n_{\mathrm{P}}\right)}{\gamma(n_{\mathrm{P}})}}
        \label{eq:xi_dot_Spyrou_lower}
    \end{equation}
    If Eq.~\ref{eq:xi_dot_Spyrou_lower} satisfies the following condition:
    \begin{equation}
        \frac{\mathrm{d}y}{\mathrm{d}t}=0 \quad \text{at} \quad \left\{\begin{aligned}
        & y=y_1\\
        & y=y_1-2 k \pi
        \end{aligned}\right.,
    \end{equation}
    then, Eq.~\ref{eq:xi_dot_Spyrou_lower} becomes the heteroclinic orbit, and such a condition is the heteroclinic bifurcation point, where $y_1$ indicates the saddle.
    \begin{equation}
        y_1 = (2k -1) \pi - \sin ^{-1} \bar{r}(n_{\mathrm{P}})
        %- \sin ^{-1} \frac{T_e\left(n_{\mathrm{P}}\right)-R}{f_{\mathrm{W}}}
    \end{equation}
    Here, $k$ is the arbitrary constant, and stable equilibrium point $y_2$ is formulated as follows (however, it has not been used in subsequent analyses):
    \begin{equation}
        y_2 = 2k \pi+\sin ^{-1} \bar{r}(n_{\mathrm{P}})
        %\frac{T_e\left(n_{\mathrm{P}}\right)-R}{f_{\mathrm{W}}}
    \end{equation}
    Finally, the condition of the heteroclinic bifurcation can be obtained as
    \begin{equation}
        \frac{1}{\bar{r}^2(n_{\mathrm{P}})} = 1 + \frac{1}{4 \gamma^2(n_{\mathrm{P}})}
    \end{equation}
    which results in
    \begin{equation}
        f_{\mathrm{W}} = \left( R-T_{\mathrm{e}}(n_{\mathrm{P}}) \right) \sqrt{\frac{k_{\mathrm{W}}^2\left(m+m_x\right)^2}{4 \gamma^2(n_{\mathrm{P}})}+1}
    \end{equation}
    On the other hand, in the case of $\frac{\mathrm{d}y}{\mathrm{d}t}>0$, then the solution becomes:
    \begin{equation}
        \frac{\mathrm{d}y}{\mathrm{d}t} = \sqrt{c_1 e^{-2 \gamma(n_{\mathrm{P}}) y}+\frac{2(\cos y-2 \gamma(n_{\mathrm{P}}) \sin y)}{1+4 \gamma^2(n_{\mathrm{P}})}+\frac{\bar{r}(n_{\mathrm{P}})}{\gamma(n_{\mathrm{P}})}}
        \label{eq:xi_dot_Spyrou_upper}
    \end{equation}
    The same approach results in the upper-speed heteroclinic bifurcation condition, which is written as
    \begin{equation}
        \begin{gathered}
            \frac{1}{\bar{r}^2(n_{\mathrm{P}})} = 1 + \frac{1}{4 \gamma^2(n_{\mathrm{P}})} \\ \text{or} \\
            \quad f_{\mathrm{W}} = \left(T_{\mathrm{e}}\left(n_{\mathrm{P}}\right)-R\right) \sqrt{\frac{k^2\left(m+m_x\right)^2}{4 \gamma^2(n_{\mathrm{P}})}+1}
        \end{gathered}
    \end{equation}
    \color{black}

\section{Method based on the cubic approximation of ``restoring'' component ~\cite{maki2014surfriding}}\label{sec:cubic_Hamiltonian}

    One of the nonlinear components that complicates the problem is the nonlinear term of the damping term. Therefore, Maki et al.~\cite{maki2014surfriding} introduced the linear approximation of this component as follows:
    \begin{equation}
        \frac{\mathrm{d}^2 y}{\mathrm{d} \tau^2}+\tilde{\beta}(n_{\mathrm{P}}) \frac{\mathrm{d} y}{\mathrm{~d} \tau}+\sin y = \bar{r}(n_{\mathrm{P}})
        \label{eq:EoM_with_liner_damping}
    \end{equation}
    where
    \begin{equation*}
        \tilde{\beta}(n_{\mathrm{P}}) = -\frac{\displaystyle \sum_{j=1}^{n_{\mathrm{M}}} \bar{A}_j(n_{\mathrm{P}}) \displaystyle \int_0^{v_\mathrm{e}} v^{j + 1} \mathrm{d}v }{\displaystyle\int_0^{v_\mathrm{e}} v^2 \mathrm{d}v} 
    \end{equation*}
    Owing to the nonlinearity of the sinusoidal function, which represents the restoring moment, Eq.~\ref{eq:EoM_with_liner_damping} is still difficult to solve analytically. Here, this sinusoidal function is approximated using a third-order polynomial, as follows:
    \begin{equation}
            \sin y \approx-\mu y\left(y-y_1\right)\left(y+y_1\right)
    \end{equation}
    where $\mu = \frac{8}{3 \pi^3}$ and $y_1=\pi$. Then, Eq.~\ref{eq:EoM_with_liner_damping} becomes
    \begin{equation}
        \frac{\mathrm{d}^2 y}{\mathrm{d} \tau^2} + \tilde{\beta}(n_{\mathrm{P}}) \frac{\mathrm{d} y}{\mathrm{d} \tau} -\mu y\left(y-y_1\right)\left(y+y_1\right) = \bar{r}(n_{\mathrm{P}})
    \end{equation}
    Now, because of this approximation, the periodicity of wave-induced force disappears; however, the approximation for one wave is sufficient for approximately obtaining the heteroclinic orbit joining two saddles. Here, assuming $a1<a2<a3$, an analytical factorization, such as Cardano's technique, yields
    \begin{equation}
        \frac{\mathrm{d}^2 y}{\mathrm{d} \tau^2} + \tilde{\beta}(n_{\mathrm{P}}) \frac{\mathrm{d} y}{\mathrm{d} \tau} -\mu\left(y-a_1\right)\left(y-a_2\right)\left(y-a_3\right)=0
        \label{eq:nonHamiltonian_factorized}
    \end{equation}
    Now, we introduce a new variable, $\tilde{y}$ as follows:
    \begin{equation}
        \frac{\mathrm{d}^2 \tilde{y}}{\mathrm{d} \tau^2} + \tilde{\beta}(n_{\mathrm{P}}) \frac{\mathrm{d} \tilde{y}}{\mathrm{d} \tau} + \tilde{\mu} \tilde{y} (1-\tilde{y}) (\tilde{y}-\tilde{a})=0
        \label{eq:EoM_for_x}
    \end{equation}
    where
    \begin{equation*}
        \left\{\begin{aligned}
                \tilde{y}&=\frac{y-a_1}{a_3-a_1}\\
                \tilde{a}&=\frac{a_2-a_1}{a_3-a_1} \\
            \tilde{\mu}&=\mu\left(a_3-a_1\right)^2
            \end{aligned}
            \right.
    \end{equation*}
    The connection between Eq.~\ref{eq:EoM_for_x} and the equation by FitzHugh-Nagumo ~\cite{nagumo1962active,mckean1970nagumo}, except for some of the coefficients, has already been discussed by Maki et al.~\cite{maki2010melnikov}. Now, the following solution form of the heteroclinic orbit is assumed~\cite{maki2010melnikov,maki2014melnikov}: 
    \begin{equation}
        \frac{\mathrm{d} \tilde{y}}{\mathrm{d} \tau} = \tilde{c} \tilde{y} (1-\tilde{y}) \quad \text{or} \quad \frac{\mathrm{d}^2 \tilde{y}}{\mathrm{d} \tau^2}=\tilde{c}^2 (1-\tilde{y}) (1-2\tilde{y})
        \label{eq:heteroclinic_orbit_cubic}
    \end{equation}
    The left side of Eq.~\ref{eq:heteroclinic_orbit_cubic} indicates that the heteroclinic orbit becomes a parabolic function on $\tilde{y}$ and the $\frac{\mathrm{d}\tilde{y}}{\mathrm{d}\tau}$ plane. Now, the substitution of Eq.~\ref{eq:heteroclinic_orbit_cubic} into Eq.~\ref{eq:EoM_for_x} yields
    \begin{equation}
        \tilde{y} \left(\tilde{\mu}-2 \tilde{c}^2\right)+\left(\tilde{c}^2+\tilde{\beta}(n_{\mathrm{P}}) \tilde{c}-\tilde{\mu} \tilde{a}\right)=0
    \end{equation}
    If the trajectory (Eq.~\ref{eq:heteroclinic_orbit_cubic}) represents the solution of Eq.~\ref{eq:EoM_for_x} for $x \in (0, 1)$, then the following condition must be satisfied:
    \begin{equation}
        \left\{\begin{array}{l}
        \tilde{\mu}-2 \tilde{c}^2=0 \\
        \tilde{c}^2+\tilde{\beta}(n_{\mathrm{P}}) \tilde{c}-\tilde{\mu} \tilde{a}=0
        \end{array}\right.
        \label{eq:cond_mu}
    \end{equation}
    The elimination of $\tilde{c}$ from Eq.~\ref{eq:cond_mu} yields
    \begin{equation}
        \frac{\tilde{\mu}}{2} \pm \tilde{\beta}(n_{\mathrm{P}}) \sqrt{\frac{\tilde{\mu}}{2}}-\tilde{\mu} \tilde{a}=0
        \label{eq:cond_myu}
    \end{equation}
    Equation~\ref{eq:cond_myu} can be solved using numerical iteration methods, such as Newton's method with respect to propeller revolution number, $n$. Then, $\tilde{c}$ or $\tilde{\mu}$ can be determined as
    \begin{equation}
        \tilde{c} = \mp \sqrt{\tilde{\mu}/2}
    \end{equation}
    Here, the heteroclinic orbit in the time domain becomes
    \begin{equation}
        \begin{gathered}
            \tilde{y}(\tau) = \frac{1}{1+\exp (-(\tilde{c} \tau -\tilde{d}))} \\
            \text{or} \\
            \tilde{y}(\tau)=\frac{1}{2}\left(1+\tanh \frac{\tilde{c} \tau -\tilde{d}}{2}\right)
        \end{gathered}
    \end{equation}
    In Eq. 59, $\tilde{d} \in (-\infty, \infty)$ is an arbitrary constant that can be determined by the initial condition.

\section{Continuous, piecewise linear approximation method~\cite{maki2010surfriding}}
    One of the nonlinearities that makes the theoretical approach difficult is the sinusoidal component included in the equation of motion (Eq.~\ref{eq:EoM_with_liner_damping}). Therefore, Maki et al.~\cite{maki2010surfriding} tried to approximate this term by using the continuous piecewise linear function, Eq.~\ref{eq:CPL_approx}:  
    %
    %\begin{equation}
    %    \begin{aligned}
    %        \sin \left(k \xi_G\right) \approx S\left(k \xi_G\right) \equiv \left\{\begin{array}{ccl}
    %        -\frac{4}{\lambda}\left(\xi_G+\frac{1}{2} \lambda\right) & : \text { Range } 1 & {\left[-\frac{3}{4} \lambda \leq \xi_G \leq-\frac{1}{4} \lambda\right]} \\
    %        \frac{4}{\lambda}\left(\xi_G+\lambda\right) & : \text { Range } 2 & {\left[-\frac{5}{4} \lambda \leq \xi_G \leq-\frac{3}{4} \lambda\right]} \\
    %        -\frac{4}{\lambda}\left(\xi_G+\frac{3}{2} \lambda\right) & : \text { Range } 3 & {\left[-\frac{7}{4} \lambda \leq \xi_G \leq-\frac{5}{4} \lambda\right]}
    %        \end{array}\right.
    %    \end{aligned}
    %    \label{eq:CPL_approx}
    %\end{equation}
    %
    \begin{equation}
        \begin{aligned}
            &\sin y \approx S\left( y \right) \\
            &\equiv \left\{\begin{array}{ccl}
            -\frac{2}{\pi}\left( y + \pi \right) & : \text { Range } 1 & {\left[-\frac{3}{2} \pi \leq y \leq -\frac{1}{2} \pi \right]} \\
            \frac{2}{\pi}\left( y + 2\pi \right) & : \text { Range } 2 & {\left[-\frac{5}{2} \pi \leq y \leq -\frac{3}{2} \pi \right]} \\
            -\frac{2}{\pi}\left( y + 3\pi \right) & : \text { Range } 3 & {\left[-\frac{7}{2} \pi \leq y \leq -\frac{5}{2} \pi \right]}
            \end{array}\right.
        \end{aligned}
        \label{eq:CPL_approx}
    \end{equation}
    \begin{figure}[htb]
        \centering
        \includegraphics[clip,width=1.0\hsize]{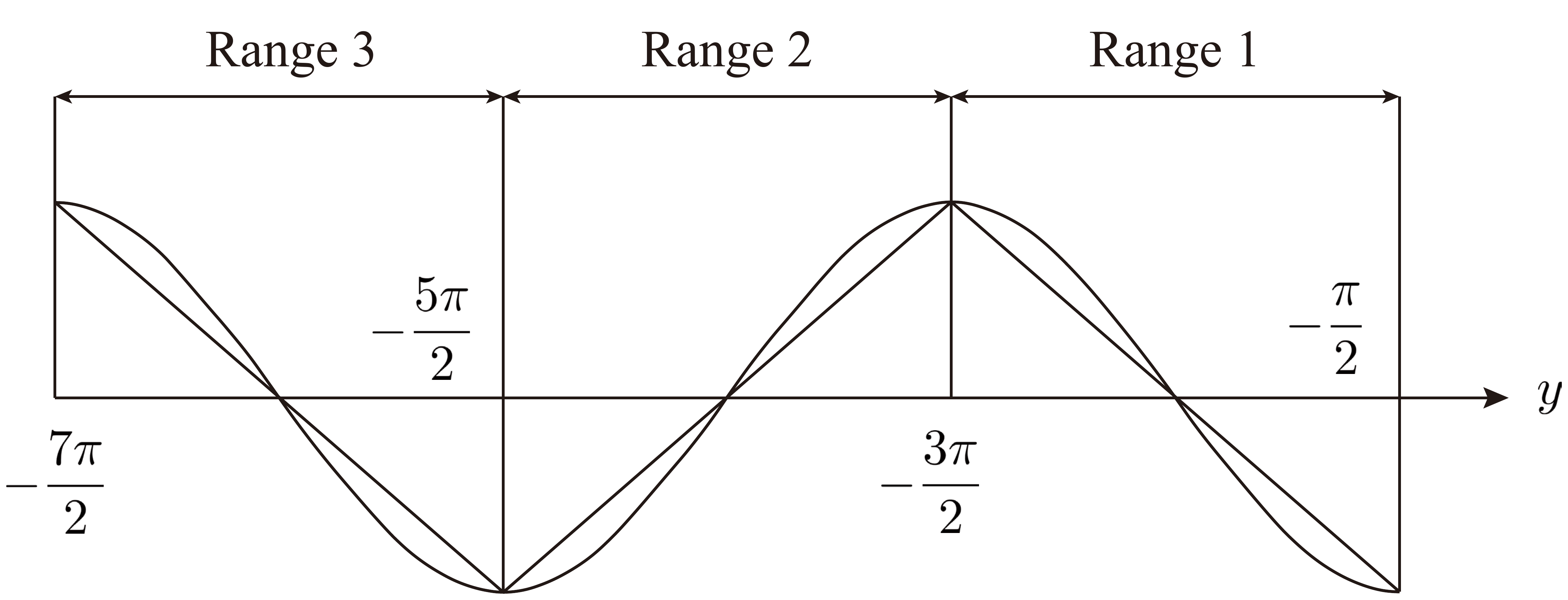}
        \caption{Schematic of piecewise linearization of a sinusoidal function}
        \label{fig:CPL_approx}
    \end{figure}
    Fig.~\ref{fig:CPL_approx} shows its schematic view. Then, Eq.~\ref{eq:EoM_with_liner_damping} becomes as follows:
    %\begin{equation}
    %    \begin{gathered}
    %        \ddot{\xi}_G+\alpha_1 \dot{\xi}_G+\alpha_2 S\left(k \xi_G\right)=\alpha_3\\
    %        \text{where} \quad \left\{ \begin{aligned}
    %            \alpha_1 &\equiv \beta(n) /\left(m+m_x\right) \\ \alpha_2 &\equiv f /\left(m+m_x\right) \\ \alpha_3 &\equiv \left[T_e(c, n)-R(c)\right] /\left(m+m_x\right)
    %        \end{aligned}\right.
    %    \end{gathered}
    %\end{equation}
    \begin{equation}
        \frac{\mathrm{d}^2 y}{\mathrm{d} \tau^2}+\tilde{\beta}(n_{\mathrm{P}}) \frac{\mathrm{d} y}{\mathrm{~d} \tau} + S(y) = \bar{r}(n_{\mathrm{P}})
        \label{eq:EoM_with_liner_damping_PL}
    \end{equation}
    Komuro~\cite{komuro1992bifurcation} as well as Endo and Chua~\cite{endo1993piecewise} have extensively analyzed the piecewise linear system. In the field of naval architecture and ocean engineering, Belenky~\cite{Belenky1993,Belenky1994} proposed a method for calculating the capsize probability in beam sea conditions based on a piecewise linear approach.
    Now, the equation of motion is linearized in each region. The solution of Range 1 is given as
    \begin{equation}
            \xi_{G 1} (\tau)=c_1 e^{\lambda_1 \tau}+c_2 e^{\lambda_2 \tau}-\frac{\pi}{2} (\bar{r}(n_{\mathrm{P}}) + 2)
        \label{eq:solution_range1}
    \end{equation}
    where
    \begin{equation*}
        \lambda_{1,2} \equiv \frac{-\alpha_1 \pm \sqrt{\alpha_1^2+16 \alpha_2 / \lambda}}{2}
    \end{equation*}
    In addition, the solution of Range 2 is formulated as
    \begin{equation}
        \begin{aligned}
            \xi_{G 2}(\tau)=c_3 e^{\lambda_3 \tau}+c_4 e^{\lambda_4 \tau} + \frac{\pi}{2} (\bar{r}(n_{\mathrm{P}}) - 4)
        \end{aligned}
        \label{eq:solution_range2}
    \end{equation}
    where
    \begin{equation*}
        \lambda_{3,4} \equiv \frac{-\alpha_1 \pm \sqrt{\alpha_1^2-16 \alpha_2 / \lambda}}{2}
    \end{equation*}
    The solution of Range 3 is formulated as
    \begin{equation}
        \begin{gathered}
            %\xi_{G 3}(t)=c_5 e^{\lambda_1 t}+c_6 e^{\lambda_2 t}-\frac{3}{2} \lambda-\frac{\lambda \alpha_3}{4 \alpha_2}
            \xi_{G 3}(\tau)=c_5 e^{\lambda_1 \tau}+c_6 e^{\lambda_2 \tau}-\frac{\pi}{2} (\bar{r}(n_{\mathrm{P}}) + 6)
        \end{gathered}
        \label{eq:solution_range3}
    \end{equation}
    Then, the solutions of Ranges 1 and 2 connect at the border as
    %\begin{equation}
    %    \left(\xi_{\mathrm{G}}, \dot{\xi}_{\mathrm{G}} \right) = \left(-\frac{3}{4} \lambda, Z_1 \right)
    %\end{equation}
    \begin{equation}
        \left( y, \frac{\mathrm{d} y}{\mathrm{d} \tau} \right) = \left(-\frac{3}{2} \pi, Z_1 \right)
        \label{eq:cond_z1}
    \end{equation}
    Here, $Z_1$ is an unknown parameter to be determined later. Then, if we impose the satisfaction of Eq.~\ref{eq:cond_z1} at time $\tau=0$ for Eq.~\ref{eq:solution_range1}, then unknown coefficients $c_1$ and $c_2$ can be determined as
    \begin{equation}
        \left\{\begin{array}{l}
        %c_1=\left(Z_1+\frac{1}{4} \lambda \lambda_2-\frac{\lambda \alpha_3}{4 \alpha_2} \lambda_2\right) /\left(\lambda_1-\lambda_2\right) \\
        %c_2=\left(-Z_1-\frac{1}{4} \lambda \lambda_1+\frac{\lambda \alpha_3}{4 \alpha_2} \lambda_1\right) /\left(\lambda_1-\lambda_2\right)
        c_1 = \frac{1}{\lambda_1-\lambda_2} \left[Z_1 + \frac{\pi}{2} \lambda_2 (1 - \bar{r}(n_{\mathrm{P}}))\right] \\
        c_2 = \frac{1}{\lambda_1-\lambda_2} \left[-Z_1 - \frac{\pi}{2} \lambda_1 (1 - \bar{r}(n_{\mathrm{P}}))\right]
        \end{array}\right.
    \end{equation}
    Similarly, the solution of Ranges 2 and 3 connect at the border as
    \begin{equation}
        \left( y, \frac{\mathrm{d} y}{\mathrm{d} t} \right) = \left(-\frac{5}{2} \pi, Z_2 \right)
        \label{eq:cond_z2}
    \end{equation}  
    Then, if we impose the satisfaction of Eq.~\ref{eq:cond_z2} at time $\tau=0$ for Eq.~\ref{eq:solution_range2}, unknown coefficients $c_5$ and $c_6$ can also be determined as
    \begin{equation}
        \left\{\begin{array}{l}
        %c_5 = \left(Z_2-\frac{1}{4} \lambda \lambda_2-\frac{\lambda \alpha_3}{4 \alpha_2} \lambda_2\right) /\left(\lambda_1-\lambda_2\right) \\
        %c_6 = \left(-Z_2+\frac{1}{4} \lambda \lambda_1+\frac{\lambda \alpha_3}{4 \alpha_2} \lambda_1\right) /\left(\lambda_1-\lambda_2\right)
        c_5 = \frac{1}{\lambda_1 - \lambda_2} \left[Z_2 - \frac{\pi}{2} \lambda_2 (1 + \bar{r}(n_{\mathrm{P}}))\right] \\
        c_6 = \frac{1}{\lambda_1 - \lambda_2} \left[-Z_2 + \frac{\pi}{2} \lambda_1 (1 + \bar{r}(n_{\mathrm{P}}))\right]
        \end{array}\right.
    \end{equation}
    Considering the heteroclinic orbit, the following condition must hold:
    \begin{equation}
        c_2 = 0 \quad \text{and} \quad c_5 = 0
    \end{equation}
    Then, unknown coefficients $Z_1$ and $Z_2$ can be determined as follows:
    \begin{equation}
        \left\{\begin{aligned}
            Z_1 &= - \frac{\pi}{2}\lambda_1 (1 - \bar{r}(n_{\mathrm{P}}))\\
            Z_2 &= \frac{\pi}{2}\lambda_2 (1 + \bar{r}(n_{\mathrm{P}}))
        \end{aligned}\right.
    \end{equation}
    Now, if we impose the satisfaction of Eq.~\ref{eq:cond_z1} at time $\tau=0$  for Eq.~\ref{eq:solution_range3}, unknown coefficients $c_1$ and $c_2$ can be determined as
    \begin{equation}
        \left\{\begin{aligned}
            c_3 &= -\frac{\pi}{2} \frac{\lambda_1+\lambda_4}{\lambda_3-\lambda_4}(1 - \bar{r}(n_{\mathrm{P}}))\\
            c_4 &= \frac{\pi}{2} \frac{\lambda_1+\lambda_3}{\lambda_3-\lambda_4}(1 - \bar{r}(n_{\mathrm{P}}))
        \end{aligned}\right.
    \end{equation}
    The matching of the trajectory after $\bar{\tau}$ seconds results in the following condition:
    \begin{equation}
        \left\{\begin{aligned}
            - \frac{\pi}{2} (1 + \bar{r}(n_{\mathrm{P}})) &= c_3 e^{\lambda_3 \bar{\tau}}+c_4 e^{\lambda_4 \bar{\tau}} \\
            \frac{\pi}{2} \lambda_2 (1 + \bar{r}(n_{\mathrm{P}})) &= c_3 \lambda_3 e^{\lambda_3 \bar{\tau}}+c_4 \lambda_4 e^{\lambda_4 \bar{\tau}}
        \end{aligned}\right.
    \end{equation}
    Here, 
    \begin{equation}
        \left\{\begin{aligned}
        & c_{3, 4} \equiv c_{\mathrm{R}} \pm i c_{\mathrm{I}} \\
        & \lambda_{3, 4} \equiv \lambda_{\mathrm{R}} \pm i \lambda_{\mathrm{I}}
        \end{aligned}\right.
    \end{equation}
    Then, the authors finally obtained the following condition, which must be satisfied by $n$ and $\bar{\tau}$; both properties are determined by Newton's method.
    \begin{equation}
        \left\{\begin{aligned}
        &-\frac{\pi}{2} (1 + \bar{r}(n_{\mathrm{P}})) \\
        &= 2 e^{\lambda_{\mathrm{R}} \bar{\tau}}\left[c_{\mathrm{R}} \cos \lambda_{\mathrm{I}} \bar{\tau} - c_{\mathrm{I}} \sin \lambda_{\mathrm{I}} \bar{\tau}\right] \\
        &\frac{\pi}{2} \lambda_2 (1 + \bar{r}(n_{\mathrm{P}})) \\
        &= 2 e^{\lambda_{\mathrm{R}} \bar{\tau}}\left[\left(c_{\mathrm{R}} \lambda_{\mathrm{R}}-c_{\mathrm{I}} \lambda_{\mathrm{I}} \right) \cos \lambda_{\mathrm{I}} \bar{\tau}-\left(c_{\mathrm{R}} \lambda_{\mathrm{I}} + c_{\mathrm{I}} \lambda_R\right) \sin \lambda_{\mathrm{I}} \bar{\tau} \right]
        \end{aligned}\right.
    \end{equation}
    
\section{Melnikov's method for the Hamiltonian part ~\cite{maki2010surfriding}}\label{sec:Melnikov_Hamiltonian_part}
    %~\cite{Kan_1990_surfriding,spyrou2006,maki2010surfriding}
    Melnikov's method ~\cite{Guckenheimer_Holmes1983,Greenspan_Holmes1984,holmes1980averaging} is a powerful analytical method to estimate the heteroclinic point of nonautonomous systems. In the field of naval architecture engineering, this method has been utilized to predict the capsizing event in seas (e.g., ~\cite{maki2010melnikov}). Now, assumed the following general nonautonomous system:
    \begin{equation}
        \frac{\mathrm{d}x}{\mathrm{d}\tau} = f(x) + \varepsilon g(x,\tau)
        \label{eq:Melnikov_explanation_nonautonoumous}
    \end{equation}
    where $x\in \mathbb{R}^n$, $f(x)$: $\mathbb{R}^n \rightarrow \mathbb{R}^n$, and $g(x, \tau)$: $\mathbb{T} \times \mathbb{R}^n \rightarrow \mathbb{R}^n$.
    Assume that the primary part of Eq~\ref{eq:Melnikov_explanation_nonautonoumous}
    \begin{equation}
        \frac{\mathrm{d}x}{\mathrm{d}\tau} = f(x) 
    \end{equation}
    has saddles. Then, we assume that the primary part has a separatrix, $x_0$, in a certain parameter condition. Then, Menikov's function is formulated as
    \begin{equation}
        M = - \int_{-\infty}^{\infty} f(x_0) \wedge \varepsilon g(x_0,\tau) \mathrm{d}\tau
    \end{equation}
    Here, $\wedge$ indicates the wedge product. In the two dimensional case, $a \wedge b = a_1 b_2 - a_2 b_1$ for $a^{\top}=(a_1,a_2) \in \mathbb{R}^2$ and $b^{\top}=(b_1,b_2) \in \mathbb{R}^2$. This indicates the projection of the difference between stable and unstable branches onto the hyperplane, which is diagonal to $x_0$.   

\begin{figure*}
        \centering
        \includegraphics[clip,width=0.85\hsize]{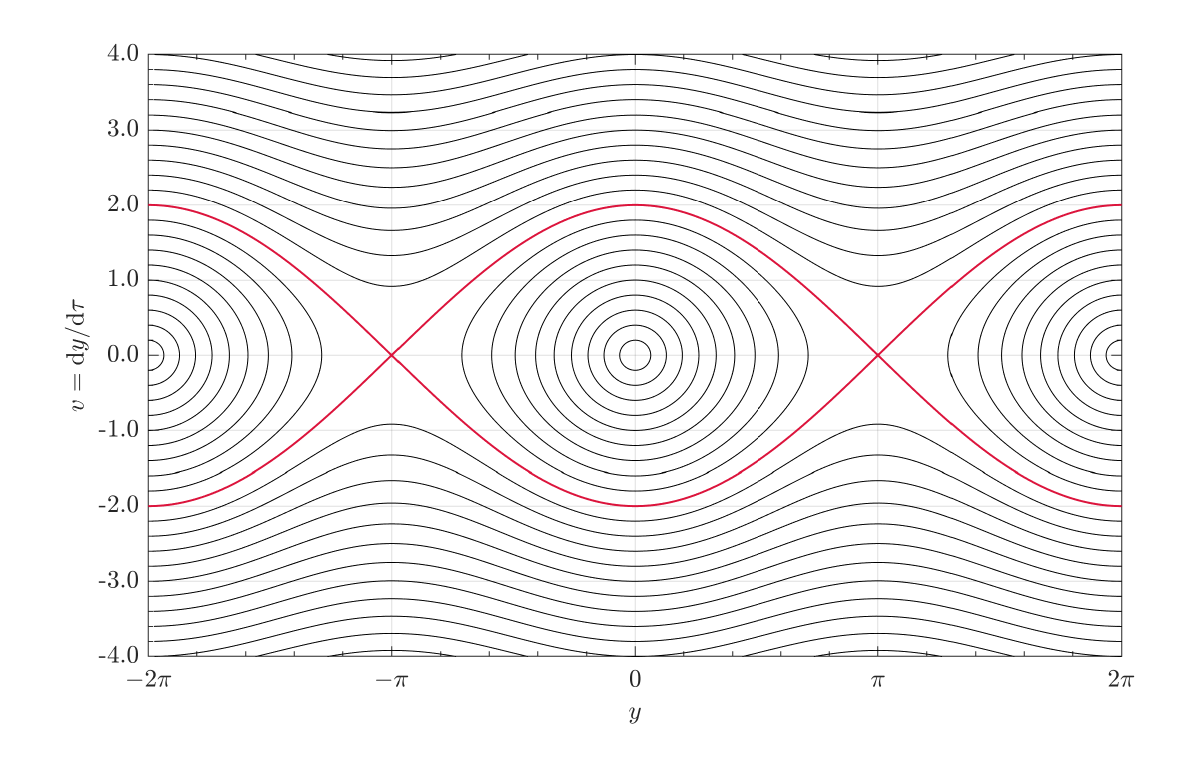}
        \caption{Phase portrait of the Hamiltonian system ($\frac{\mathrm{d}^2 y}{\mathrm{d} \tau^2} + \sin y = 0$). Here, the red solid lines indicate the separatrix (Eq.~\ref{eq:separatrix_for_Melnikov}).}
        \label{fig:phase_portrait_Hamiltonian_system}
    \end{figure*}

    Now, let $x \in \mathbb{R}^2$ be
    \begin{equation}
        x(\tau) \equiv (y, \frac{\mathrm{d}y}{\mathrm{d}\tau})^{\top}
    \end{equation}
    Then, we define $f(x)\in \mathbb{R}^2 \rightarrow \mathbb{R}^2$ and $\varepsilon g(x,\tau)\in \mathbb{R}^2 \rightarrow \mathbb{R}^2$ for the components of Eq.~\ref{eq:equation_of_motion_nondim_2} as follows:
    \begin{equation}
        \frac{\mathrm{d}x}{\mathrm{d}\tau} = f(x) + \varepsilon g(x,\tau)
    \end{equation}
    where
    \begin{equation*}
        \left\{\begin{aligned}
            f(x) &\equiv \left(\frac{\mathrm{d} y}{\mathrm{d} \tau}, - \sin y \right)^{\top}\\
            g(x) &\equiv \left(0,  \bar{r}\left(n_{\mathrm{P}}\right) - \sum_{k=1}^{n_{\mathrm{M}}} \bar{A}_k(n) \left( \frac{\mathrm{d} y}{\mathrm{d} \tau} \right)^k \right)^{\top} 
        \end{aligned}\right.
    \end{equation*}
    The Hamiltonian part of Eq.~\ref{eq:equation_of_motion_nondim_2} is defined as
    \begin{equation}
        \frac{\mathrm{d}^2 y}{\mathrm{d} \tau^2} + \sin y = 0
        \label{eq:Hamiltonian_system}
    \end{equation}
    The phase portrait of the Hamiltonian system (Eq.~\ref{eq:Hamiltonian_system}) is shown in Fig. \ref{fig:phase_portrait_Hamiltonian_system}. 
    The trajectory connecting $y = \pm \pi$ on the lower or upper sides of vector field can be determined as follows:
    \begin{equation}
        v \equiv \frac{\mathrm{d} y}{\mathrm{d} \tau} = \mp 2 \cos \frac{y}{2}.
        \label{eq:separatrix_for_Melnikov}
    \end{equation}
    This is called a ``separatrix.'' 
    Here, the Melnikov function is defined as follows:
    \begin{equation}
        M \equiv \int_{\mp \infty}^{\pm \infty} v \left( 
        \bar{r}\left(n_{\mathrm{P}}\right) - \sum_{k=1}^{n_{\mathrm{M}}} \bar{A}_k(n) \left( \frac{\mathrm{d} y}{\mathrm{d} \tau} \right)^k
        \right) \mathrm{d} \tau
        \label{eq:Malenikov_original_cal}
    \end{equation}
    By substituting Eq.~\ref{eq:separatrix_for_Melnikov} into Eq.~\ref{eq:Malenikov_original_cal}, we get
    \begin{equation}
        M = \int_{\mp \pi}^{\pm \pi} \bar{r}\left(n_{\mathrm{P}}\right) \mathrm{~d} y - \sum_{k=1}^{n_{\mathrm{M}}} \bar{A}_k(n) \int_{\mp \pi}^{\pm \pi} \left( \mp 2 \cos \frac{y}{2} \right)^k \mathrm{~d} y
    \end{equation}
    Assuming that $M=0$, the following relationship can be obtained:
    \begin{equation}
        2 \pi \bar{r}\left(n_{\mathrm{P}}\right) = \sum_{k=1}^{n_{\mathrm{M}}} \bar{A}_k(n) \left( \mp 2 \right)^k I_k
    \end{equation}
    where
    \begin{equation*}
        I_k \equiv \int_{-\pi}^\pi \cos ^k(y / 2) \mathrm{d} y
    \end{equation*}
    Here, $I_k$ can be evaluated using the Gamma function, $\Gamma$.
    \begin{equation}
        I_k = 2 \sqrt{\pi} \Gamma\left(\frac{k+1}{2}\right) / \Gamma\left(\frac{k+2}{2}\right)
    \end{equation}
    Note that $I_1 = 4$, $I_2 = \pi$, $I_3 = \frac{8}{3}$, $I_4 = \frac{3}{4}\pi$, and $I_5 = \frac{32}{15}$. Here, if we put $n_{\mathrm{M}}=1$, then we have
    \begin{equation}
        \begin{gathered}
        \frac{T_e\left(c_{\mathrm{W}} ; n_{\mathrm{P}}\right)-R\left(c_{\mathrm{W}}\right)}{f} = - \frac{4 c_1}{\pi \sqrt{f_{\mathrm{W}} k_{\mathrm{W}} (m+m_x)}}
        \end{gathered}
    \end{equation}
    This relation is completely identical to the formula obtained by Kan ~\cite{Kan_1990_surfriding}.
    Next, if we set $n_{\mathrm{M}}=3$ and $\kappa_3 = 0$, then we have
    \begin{equation}
        \begin{aligned}
        &\frac{T_e\left(c_{\mathrm{W}} ; n_{\mathrm{P}}\right)-R\left(c_{\mathrm{W}}\right)}{f} = - \frac{4 \left(c_1 + 2 c_2 c_{\mathrm{W}} + 3c_3 c_{\mathrm{W}}^2 \right)}{\pi \sqrt{f_{\mathrm{W}} k_{\mathrm{W}} (m+m_x)}}\\
        &+ \frac{2 \left(c_2 + 3c_3 c_{\mathrm{W}} \right)}{k_{\mathrm{W}} (m+m_x)} - \frac{32 c_3 \sqrt{f_{\mathrm{W}}}}{3 \pi \left[ k_{\mathrm{W}} (m+m_x) \right]^{\frac{3}{2}}}
        \end{aligned}
    \end{equation}
    This relation is completely identical to the formula obtained by Spyrou \cite{spyrou2006}.
    
    Finally, if we set $n_{\mathrm{M}}=5$ and $\kappa_j = 0~(j=3 \sim 5)$, then we have
    \begin{equation}
        \begin{aligned}
        &\frac{T_e\left(c_{\mathrm{W}} ; n_{\mathrm{P}}\right)-R\left(c_{\mathrm{W}}\right)}{f} \\
        = &- \frac{4 \left(c_1 + 2 c_2 c_{\mathrm{W}} + 3c_3 c_{\mathrm{W}}^2 + 4c_4 c_{\mathrm{W}}^3 + 5c_5 c_{\mathrm{W}}^4 \right)}{\pi f_{\mathrm{W}}^{1/2} k_{\mathrm{W}}^{3/2} (m+m_x)^{3/2}}\\
        &+ \frac{2 \left(c_2 + 3c_3 c_{\mathrm{W}} + 6c_4 c_{\mathrm{W}}^3  + 10c_5 c_{\mathrm{W}}^3 \right)}{k_{\mathrm{W}} (m+m_x)}\\
        &- \frac{32 \left( c_3 + 4c_4 c_{\mathrm{W}} + 10c_5 c_{\mathrm{W}}^2 \right) f_{\mathrm{W}}^{1/2}}{3 \pi k_{\mathrm{W}}^{3/2} (m+m_x)^{3/2}}\\
        &+ \frac{6 \left(c_4 + 5c_5 c_{\mathrm{W}} \right) f_{\mathrm{W}} }{k_{\mathrm{W}}^2 (m+m_x)^2}\\
        &- \frac{512 c_5 f_{\mathrm{W}}^{3/2}}{15 \pi k_{\mathrm{W}}^{5/2} (m+m_x)^{5/2}}
        \end{aligned}
        \label{eq:Melnikov_5th}
    \end{equation}
    This relation is completely identical to the second-generation intact stability criteria (SGISC)~\cite{Imo2020-bp,Imo2022-uf} described in Section~\ref{sec:IMO_criteria}.

\section{Melnikov's method for the non-Hamiltonian part (I) ~\cite{maki2016surfriding}}
    The original Melnikov's method is applicable for the separatrix of the Hamiltonian part of the system. In addition, the extended Melnikov's method ~\cite{salam1987mel} can be applied to the heteroclinic orbit of the non-Hamiltonian part of the system. In the field of naval architects and ocean engineering, this extended Melnikov's method is used in the study of the capsizing phenomena in regular beas seas~\cite{wu_mccue2008application,maki2010melnikov,maki2014melnikov,maki2022NOLTA}. Maki and Miyauchi~\cite{maki2016surfriding} applied the extended Melnikov's method to the heteroclinic orbits, which had been already obtained in ~\cite{spyrou2006}. This section presents the results of this method.  
    %\begin{equation}
    %    \begin{gathered}
    %        \ddot{y}+p \dot{y}^2 \operatorname{sgn} \dot{y}+l \dot{y}+m \dot{y}^2+\sin y=\bar{r(n)}\\
    %        \text{where} \quad \left\{\begin{aligned}
    %            l=\frac{a}{\sqrt{q}\left(m+m_x\right)}, \quad m=\frac{b}{k\left(m+m_x\right)}
    %        \end{aligned}\right.
    %    \end{gathered}
    %\end{equation}

    Here, Eq.~\ref{eq:equation_of_motion_spyrou} is solved using the following additional correction term:
    \begin{equation}
        \begin{aligned}
            &\frac{\mathrm{d}^2 y}{\mathrm{d} \tau^2} + \gamma(n_{\mathrm{P}}) \left( \frac{\mathrm{d} y}{\mathrm{d} \tau}  \right)^2 \operatorname{sgn} \left( \frac{\mathrm{d} y}{\mathrm{d} \tau} \right)\\
            &+ \Gamma_1(n_{\mathrm{P}}) \frac{\mathrm{d} y}{\mathrm{d} \tau} + \Gamma_2(n_{\mathrm{P}}) \left( \frac{\mathrm{d} y}{\mathrm{d} \tau} \right)^2 +\sin y = \bar{r}(n_{\mathrm{P}})
        \end{aligned}
    \end{equation}
    In the above-mentioned equation, $\Gamma_1(n_{\mathrm{P}})$ and $\Gamma_2(n_{\mathrm{P}})$ are determined as the best fit of the original nonlinear damping term in terms of the least square fit. Now, to apply the extended Melnikov's method, $\sigma$ is added to both sides to get
    \begin{equation}
        \begin{aligned}
            &\frac{\mathrm{d}^2 y}{\mathrm{d} \tau^2} + \gamma(n_{\mathrm{P}}) \left( \frac{\mathrm{d} y}{\mathrm{d} \tau}  \right)^2 \operatorname{sgn} \left( \frac{\mathrm{d} y}{\mathrm{d} \tau} \right) + \sin y - \bar{r}(n_{\mathrm{P}})+\sigma\\
            &=-\left[ \Gamma_1(n_{\mathrm{P}}) \frac{\mathrm{d} y}{\mathrm{d} \tau} + \Gamma_2(n_{\mathrm{P}}) \left( \frac{\mathrm{d} y}{\mathrm{d} \tau} \right)^2 \right]+\sigma
        \end{aligned}
    \end{equation}
    The left part of the above-mentioned equation has been explained in Section~\ref{sec:Spyrou}. If $\sigma$ satisfies the following condition, the heteroclinic trajectory is achieved:
    \begin{equation}
        \sigma = \bar{r}(n_{\mathrm{P}}) \pm \frac{2 \gamma(n_{\mathrm{P}})}{\sqrt{1+4 \gamma^2(n_{\mathrm{P}})}}
    \end{equation}
    The trajectory in the phase plane can be expressed as 
    \begin{equation}
        \begin{aligned}
             &v = \mp \frac{2}{\left(1+4 \gamma^2(n_{\mathrm{P}})\right)^{1 / 4}}\left|\cos \frac{y + \varepsilon_{y}}{2}\right| \\
             &\text { where } \tan \varepsilon_{y}=\mp 2 \gamma(n_{\mathrm{P}})
        \end{aligned}
        \label{eq:sol_Spyrou_nonHamilton}
    \end{equation}
    The extended Melnikov's integral can be evaluated as follows:
    \begin{equation}
        \begin{aligned}
            M &= -\int_{-\infty}^{\infty} v \left(\sigma-\Gamma_1(n_{\mathrm{P}}) v -\Gamma_2(n_{\mathrm{P}}) v^2 \right) e^{\mp 2 \gamma(n_{\mathrm{P}}) y} \mathrm{~d} \tau\\
              &= -\int_{-\pi-\delta}^{\pi-\delta}\left(\sigma-\Gamma_1(n_{\mathrm{P}}) v-\Gamma_2(n_{\mathrm{P}}) v^2\right) e^{\mp 2 \gamma(n_{\mathrm{P}}) y} \mathrm{~d} y
        \end{aligned}
    \end{equation}
    By imposing $M=0$ and considering Eq.\ref{eq:sol_Spyrou_nonHamilton}, the following condition of the surf-riding threshold is achieved:
    \begin{equation}
        \begin{aligned}
            &\mp \frac{\sigma\left(-1+e^{\mp 4 \gamma(n_{\mathrm{P}}) y}\right)}{\gamma(n_{\mathrm{P}})}=\\
            & \frac{8 \Gamma_1(n_{\mathrm{P}})\left(1+e^{\mp 4 \gamma(n_{\mathrm{P}}) y}\right)}{\left(1+16 \gamma^2(n_{\mathrm{P}})\right)\left(1+4 \gamma^2(n_{\mathrm{P}})\right)^{1 / 4}} \\
            &\mp \frac{2 \Gamma_2(n_{\mathrm{P}}) \left(-1+e^{\mp 4 \gamma(n_{\mathrm{P}}) y}\right)}{\gamma(n_{\mathrm{P}})\left(1+4 \gamma^2(n_{\mathrm{P}})\right) \sqrt{1+4 \gamma^2(n_{\mathrm{P}})}}
        \end{aligned}
    \end{equation}

\section{Melnikov's method for non-Hamiltonian part II~\cite{maki2016surfriding}}
\begin{figure}[h]
        \centering
        \includegraphics[clip,width=1.0\hsize]{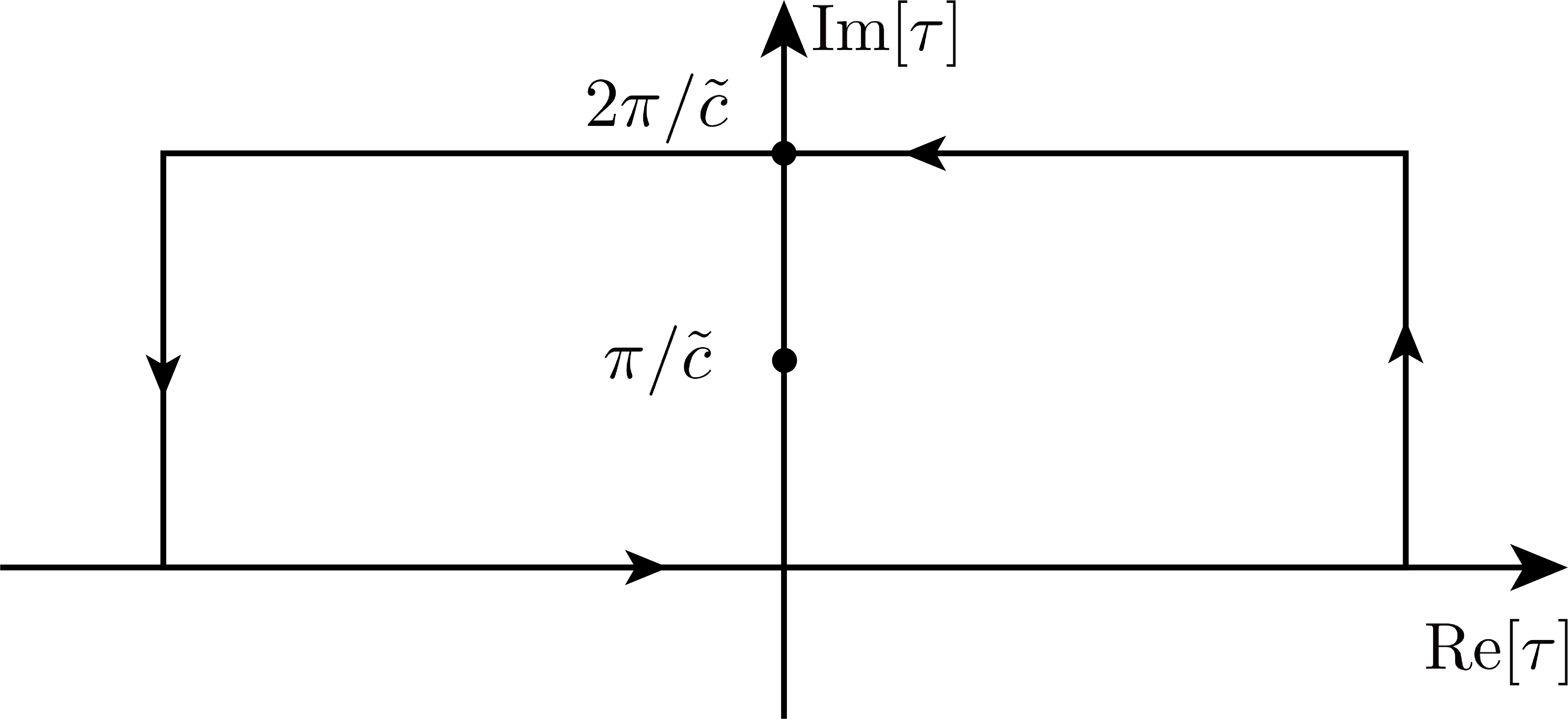}
        \caption{Integral path for calculating $I_1$. Here, $\tau=\pi i / \tilde{c}$ represents the pole of the integrand.}
        \label{fig:integral_path}
    \end{figure}

    Maki and Miyauchi also applied the extended Melnikov's method~\cite{salam1987mel} to the system used in ~\cite{maki2010surfriding}. They introduced the correction components for the approximated wave force part as having the form of a fifth-order polynomial as follows:
    \begin{equation}
        \sin y \approx-\mu y\left(y-y_1\right)\left(y-y_2\right) + \sigma_1 y + \sigma_3 y^3 + \sigma_5 y^5
    \end{equation}
    Moreover, the correction components for the damping term, ($\tilde{\beta}_1(n_{\mathrm{P}})$ and $\tilde{\beta}_2(n_{\mathrm{P}})$), were also introduced: 
    \begin{equation}
        \begin{aligned}
            \frac{\mathrm{d}^2 y}{\mathrm{d} \tau^2} + & \tilde{\beta}(n_{\mathrm{P}}) \frac{\mathrm{d} y}{\mathrm{d} \tau} + \tilde{\beta}_1(n_{\mathrm{P}}) \frac{\mathrm{d} y}{\mathrm{d} \tau} + \tilde{\beta}_2(n_{\mathrm{P}}) \left( \frac{\mathrm{d} y}{\mathrm{d} \tau} \right)^2\\
            &-\mu y\left(y-y_1\right)\left(y-y_2\right) + \sigma_1 y + \sigma_3 y^3 + \sigma_5 y^5 = \bar{r}
        \end{aligned}
    \end{equation}  
    Then, $\sigma$ is added to both sides to get
    \begin{equation}
        \begin{aligned}
            &\frac{\mathrm{d}^2 y}{\mathrm{d} \tau^2}  +\tilde{\beta}(n_{\mathrm{P}}) \frac{\mathrm{d} y}{\mathrm{d} \tau} - \mu y\left(y-y_1\right)\left(y-y_2\right)- \bar{r}(n_{\mathrm{P}}) +\sigma \\
            & =-\Bigg(\tilde{\beta}_1(n_{\mathrm{P}}) \frac{\mathrm{d} y}{\mathrm{d} \tau} + \tilde{\beta}_2(n_{\mathrm{P}}) \left( \frac{\mathrm{d} y}{\mathrm{d} \tau} \right)^2 \\
            &+ \sigma_1 y + \sigma_3 y^3 + \sigma_5 y^5\Bigg)+\sigma
        \end{aligned}
        \label{eq:EoM_nonhamiltonoan_II}
    \end{equation}
    First, the authors solve the follwoing non-Hamiltonian part of Eq.~\ref{eq:EoM_nonhamiltonoan_II}:
    \begin{equation}
        \begin{aligned}
            \frac{\mathrm{d}^2 y}{\mathrm{d} \tau^2} & +\tilde{\beta}(n_{\mathrm{P}}) \frac{\mathrm{d} y}{\mathrm{d} \tau} - \mu y\left(y-y_1\right)\left(y-y_2\right)- \bar{r}(n_{\mathrm{P}}) +\sigma = 0
        \end{aligned}
        \label{eq:EoM_nonhamiltonoan_II_mainpart}
    \end{equation}
    As shown in Section~\ref{sec:cubic_Hamiltonian}, the solution of Eq.~\ref{eq:EoM_nonhamiltonoan_II_mainpart} is written as
    \begin{equation}
        \left\{\begin{aligned}
            \tilde{y} &= \frac{1}{1+\exp (-\tilde{c} \tau + d)} \\
            \frac{\mathrm{d} \tilde{y}}{\mathrm{d} \tau} &= \frac{\exp (-\tilde{c} \tau + d)}{\left[1+\exp (-\tilde{c} \tau + d)\right]^2}
        \end{aligned}\right.
        \label{eq:sol_nonHamiltonian_part}
    \end{equation}
    where
    \begin{equation*}
        \begin{aligned}
            \tilde{y} &= \frac{y-a_1}{a_3-a_1} \quad \text{and} \quad \frac{\mathrm{d} \tilde{y}}{\mathrm{d} \tau} = \frac{1}{a_3-a_1} v
        \end{aligned}
    \end{equation*}
    Then, the Melnikov's integral can be evaluated as follows:
    \begin{equation}
        \begin{aligned}
            M=&-\int_{-\infty}^{\infty} v \Big[\sigma-\left(\tilde{\beta}_1(n_{\mathrm{P}}) v +\tilde{\beta}_2(n_{\mathrm{P}}) v^2 \right)\\
            &-\left( \sigma_1 y + \sigma_3 y^3 + \sigma_5 y^5 \right)\Big] e^{\tilde{\beta} \tau} \mathrm{~d} \tau
        \end{aligned}
    \end{equation}
    By considering Eq.~\ref{eq:sol_nonHamiltonian_part}, $M$ can be evaluated as follows:
    \begin{equation}
        \begin{aligned}
            M= &\sigma\left(a_3-a_1\right) I_1-\tilde{\beta}_1(n_{\mathrm{P}}) \left(a_3-a_1\right)^2 I_2\\
            &-\tilde{\beta}_2(n_{\mathrm{P}}) \left(a_3-a_1\right)^3 I_3 \\
            & -\left(a_3-a_1\right) \cdot\\
            &\left(a_0 K_0+a_1 K_1+a_2 K_2+a_3 K_3+a_4 K_4+a_5 K_5\right)
        \end{aligned}
    \end{equation}  
    Here, the definitions of $a_1$, $a_2$, and $a_3$ are shown in Eq.~\ref{eq:nonHamiltonian_factorized}.
    Now, integral $I_i$ can be calculated as follows, and $I_n$ is defined as
    \begin{equation}
        I_n=\tilde{c}^n \int_{-\infty}^{\infty} \frac{\exp (-n\tilde{c} \tau+\tilde{\beta}(n_{\mathrm{P}}) \tau)}{[1+\exp (-\tilde{c} \tau)]^{2n}} \mathrm{~d} \tau
    \end{equation}
    In addition, $I_1$ can be calculated for the integral path on complex domain, as shown in Fig.~\ref{fig:integral_path}

    \begin{equation}
        \begin{aligned}
            I_1&=\tilde{c} \int_{-\infty}^{\infty} \frac{\exp (-\tilde{c} \tau+\tilde{\beta}(n_{\mathrm{P}}) \tau)}{[1+\exp (-\tilde{c} \tau)]^2} \mathrm{~d} \tau\\
            &=\frac{\pi \tilde{\beta}(n_{\mathrm{P}}) \csc (\pi \tilde{\beta}(n_{\mathrm{P}}) / \tilde{c}) \operatorname{sgn} \tilde{c}}{\tilde{c}}
        \end{aligned}
    \end{equation}
    The application of integration by part yields the following asymptotic relation:
    \begin{equation}
        I_{n+1} = - \frac{\tilde{\beta}^2(n_{\mathrm{P}}) - n^2\tilde{c}^2}{2n(2n+1)\tilde{c}}I_n
    \end{equation}
    Therefore, $I_{n}$ for $n > 1$ can be evaluated as
    \begin{equation}
        \begin{aligned}
            I_n = &(-1)^{n-1} \frac{\pi \tilde{\beta}(n_{\mathrm{P}}) \csc (\pi \tilde{\beta}(n_{\mathrm{P}}) / \tilde{c}) \operatorname{sgn} \tilde{c} }{(2n-1)! \tilde{c}^2} \cdot\\
            &\displaystyle\prod_{j=2}^{n} \left( (j-1) \tilde{\beta}^2(n_{\mathrm{P}}) - \tilde{c}^2\right)
        \end{aligned}
    \end{equation}
    
    %\begin{equation}
    %    \begin{aligned}
    %    I_2 & =\tilde{c}^2 \int_{-\infty}^{\infty} \frac{\exp (-2\tilde{c} \tau+\tilde{\beta}(n_{\mathrm{P}}) \tau)}{[1+\exp (-\tilde{c} \tau)]^4} \mathrm{~d} \tau =-\frac{\pi \tilde{\beta}\left(\tilde{\beta}^2-\tilde{c}^2\right) \csc (\pi \tilde{\beta} / \tilde{c}) \operatorname{sgn} \tilde{c}}{6 \tilde{c}^2}
    %    \end{aligned}
    %\end{equation}

    %\begin{equation}
    %    \begin{aligned}
    %        I_3 & =\tilde{c}^3 \int_{-\infty}^{\infty} \frac{\exp (-3\tilde{c} \tau+\tilde{\beta}(n_{\mathrm{P}}) \tau)}{[1+\exp (-\tilde{c} \tau)]^6} \mathrm{~d} \tau = \frac{\pi \tilde{\beta}\left(\tilde{\beta}^2-4 \tilde{c}^2\right)\left(\tilde{\beta}^2-\tilde{c}^2\right) \csc (\pi \tilde{\beta} / \tilde{c}) \operatorname{sgn} \tilde{c}}{120 \tilde{c}^3}
    %    \end{aligned}
    %\end{equation}  
    %
    Now, the result of integral $K_i$ is calculated. Here, $K_n$ is defined as follows:
    \begin{equation}
        K_n =\tilde{c} \int_{-\infty}^{\infty} \frac{\exp (-\tilde{c} \tau+\tilde{\beta}(n_{\mathrm{P}}) \tau)}{[1+\exp (-\tilde{c} \tau)]^{n+2}} \mathrm{~d} \tau
    \end{equation}
    Now, $K_0=I_1$. The application of integration by part yields the following asymptotic relation:
    \begin{equation}
        K_{n+1} = \frac{(n+1)\tilde{c}+\tilde{\beta}(n_{\mathrm{P}})}{(n+2)\tilde{c}}K_n
    \end{equation}
    Therefore, $K_n$ for $n \ge 1$ can be evaluated as
    \begin{equation}
        K_n = \frac{\pi \tilde{\beta}(n_{\mathrm{P}}) \csc (\pi \tilde{\beta}(n_{\mathrm{P}}) / \tilde{c}) \operatorname{sgn} \tilde{c} }{(n+1)! \tilde{c}^{n+1}} \displaystyle\prod_{j=1}^{n} \left( j \tilde{c} + \tilde{\beta}(n_{\mathrm{P}}) \right)
    \end{equation}

\section{Comparison of the surf-riding threshold between theoretical and free-running model experiments~\cite{maki2016surfriding}}
\label{sec:comparison}
\begin{figure}[h]
        \centering
        \includegraphics[clip,width=1.0\hsize]{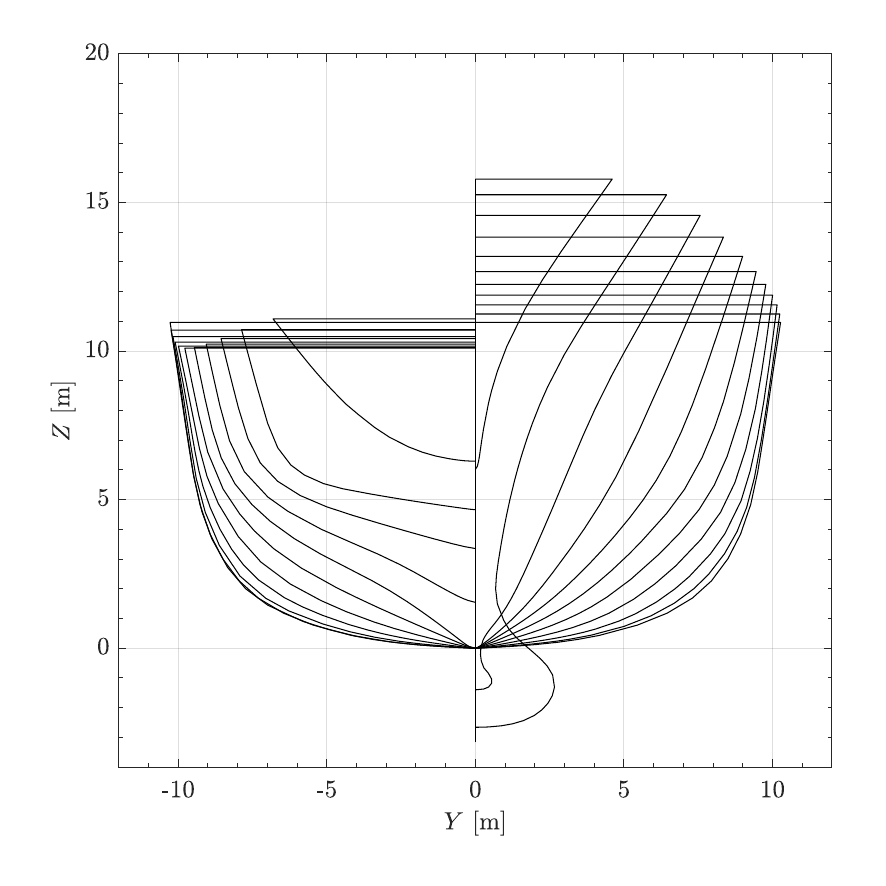}
        \caption{Body plan of DTMB5415 in real scale ($\mathrm{Length}: 142.17~\mathrm{m}$)}
        \label{fig:body_plan}
    \end{figure}

    Here, we provide the comparison of the surf-riding threshold between the theoretical and experimental free-running models. The target ship is the hull form known as DTMB5415; its body plan in real scale is shown in Fig.~\ref{fig:body_plan}. In addition, its principal particulars and coefficient list are shown in Tables~\ref{tab:principal_5415} and ~\ref{tab:coef_5415}, respectively. The experiment was conducted in the towing tank ($\text{Length}:257~\mathrm{m} \times \text{Width}:12.5~\mathrm{m}  \times \text{Depth}:7~\mathrm{m}$) of the Naval Systems Research Center (Acquisition, Technology and Logistics Agency, MINISTRY OF DEFENSE, Japan)). The free-running model is a twin-screw and twin-rudder vessel and is equipped with two propulsion motors ($200~\mathrm{W}$ output in each) and two steering motors (steering speed of $30~\mathrm{deg./s}$). Each directional angle and angular velocity were measured using a fiber-optic gyro. Based on the observed state, the model was operated in a straight line in following waves by proportional and differential (PD) control. This free-running model test is unique in that it was conducted in a towing tank. Unlike the experiment conducted in the rectangular tank, the model was able to travel a longer distance, and thus the surf-riding phenomenon could be tested up to very high speeds of approximately $\overline{\mathrm{Fn}}=0.8$ or more.

\begin{figure*}
        \centering
        \includegraphics[clip,width=1.0\hsize]{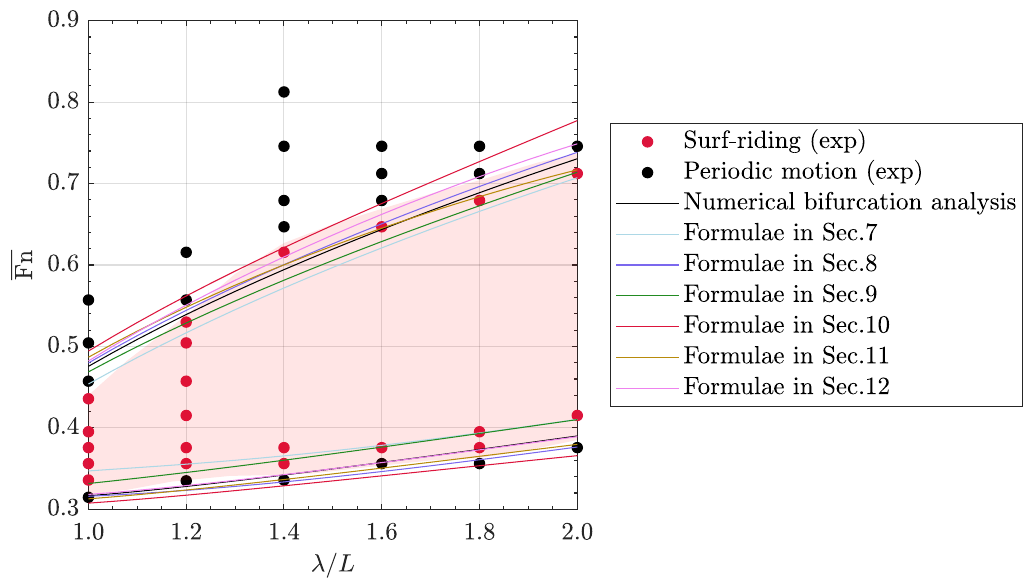}
        \caption{Comparison of the surf-riding and wave-blocking threshold for the proposed methods with the experimental results as a function of $\lambda/L$, with $H/\lambda=0.04$ (This figure duplicates Figs.~9 and 10 in the literature~\cite{maki2016surfriding}, but with minor modifications).}
        \label{fig:final_result}
    \end{figure*}

    \begin{table}[h]
        \caption{Principal particulars of DTMB5415 in $1/51.7$ scale model ($\mathrm{Length}: 2.75~\mathrm{m}$)}
        \centering
        \begin{tabular}{lc|lc} 
        item & value & item & value \\
        \hline $L_{\mathrm{WL}}~\left[ \mathrm{~m~}\right]$ & $2.750$ &  $\nabla~\left[ \mathrm{~m^3~} \right]$ & $0.0626$ \\ $B_{\mathrm{S}}~\left[ \mathrm{~m~} \right]$ & $0.369$ & $C_{\mathrm{b}}$ & $0.507$ \\
        $d_{\mathrm{S}}~\left[ \mathrm{~m~} \right]$ & $0.119$ & $C_{\mathrm{p}}$ & $0.618$ \\
        $D_{\mathrm{P}}~\left[ \mathrm{~m~} \right]$ & $0.1045$ & scale & $1/51.7$
        \end{tabular}
        \label{tab:principal_5415}
    \end{table}
    \begin{table}[h]
        \caption{Coefficient list of DTMB5415 in $1/51.7$ scale model ($\mathrm{Length}: 2.75~\mathrm{m}$)}
        \centering
        \begin{tabular}{lc|lc} 
        item & value & item & value \\
        \hline $r_1~\left[ \mathrm{~N \cdot s / m}\right]$ & $9.407$ & $\kappa_0$ & $0.6882$ \\
        $r_2~\left[ \mathrm{~N \cdot s^2 / m^2~} \right]$ & $-21.96$ & $\kappa_1$ & $-0.4047$ \\
        $r_3~\left[ \mathrm{~N \cdot s^3 / m^3~} \right]$ & $19.56$ & $\kappa_2$ & $-0.09504$ \\
        $r_4~\left[ \mathrm{~N \cdot s^4 / m^4~} \right]$ & $-5.243$ & $1-w_{\mathrm{P}}$ & $0.94$ \\
        $r_5~\left[ \mathrm{~N \cdot s^5 / m^5~} \right]$ & $0.4599$ & $1-t_{\mathrm{P}}$ & $0.85$
        \end{tabular}
        \label{tab:coef_5415}
    \end{table}
    \begin{figure}
        \centering
        \includegraphics[clip,width=1.0\hsize]{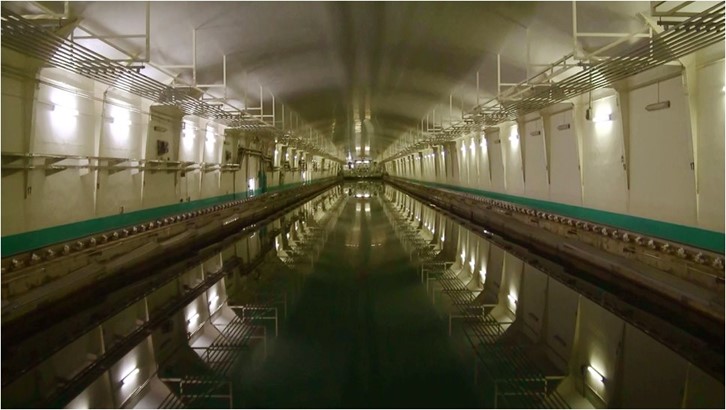}
        \caption{Towing tank of NSRC}
    \end{figure}
    \begin{figure}
        \centering
        \includegraphics[clip,width=1.0\hsize]{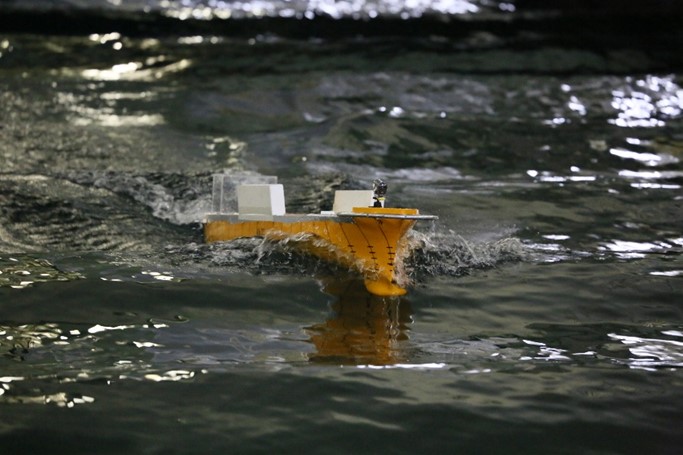}
        \caption{Free-running model experiment}
    \end{figure}
    The results based on all the theoretical and numerical bifurcation analyses described in this paper are compared with the experimentally obtained surf-riding threshold in Fig.~\ref{fig:final_result}. First, the surf-riding thresholds obtained by numerical bifurcation analysis show a good correlation with the experimentally obtained surf-riding threshold. The results obtained by numerical bifurcation analysis are the correct solutions of the surf-riding threshold of the equation of motion which we are dealing with. The equations used in this study, that is Eq.~\ref{eq:equation_of_motion}, are quite simple, and it considers only one degree-of-freedom surge motion. However, it can be seen that the equations can be used to estimate the nonlinear surge motion of the ship in the following seas with quantitative accuracy.
    
    From this figure, the trend of the theoretically obtained surf-riding and wave-blocking thresholds is generally close. However, there is some variation in the accuracy of the quantitative estimates. Among them, however, the theoretical estimation methods based on Melnikov's method are found to be able to estimate surf-riding and wave-blocking thresholds with particularly high accuracy.

\section{ SGISC}\label{sec:IMO_criteria}
    The IMO developed the SGISC to prevent ship accidents due to roll motion.
    Here, the ``intact stability'' indicates the stability of a nondamaged ship, and the antonym is ``damage stability''.
    The current draft interim guidelines on the SGISC were finalized by the IMO's Sub-Committee on Ship Design and Construction in 2020~\cite{Imo2020-bp,Imo2022-uf} and are under a trial period.
    
    The SGISC represent physics-based criteria for assessing five stability failure modes, namely, broaching, parametric roll, pure loss of stability, dead ship condition, and excessive acceleration, and possess a multilayered structure for each stability failure mode.
    The first two layers comprise level-1 and level-2 vulnerability criteria, which are simplified assessments of the failure modes.
    The third layer possesses direct stability assessment, which is a probabilistic assessment based on ship-motion equations validated with a model experiment.
    
    The broaching phenomenon is a stability failure mode that has been considered in the SGISC.
    This is because the centrifugal force due to broaching induces ship-roll motion, usually in the direction away from the center of the yaw motion.
    As surf-riding is a precursor to broaching, the current draft Interim Guidelines on the SGISC adopt level-1 and level-2 vulnerability criteria for broaching based on the surf-riding dynamics.
    This section explains the level-1 and level-2 vulnerability criteria for broaching.

    The basic concept of the vulnerability criteria shows that a ship is considered vulnerable to the broaching failure mode
    if the service speed is higher than the threshold of the nominal speed of the ship,
    at which surf-riding occurs, regardless of the initial ship speed (the critical nominal speed of the ship).
    The critical nominal speed of the ship is calculated as the nominal speed of the ship at which the heteroclinic bifurcation occurs (Sec.~\ref{section_PhasePortrait}).
    Here, the nominal speed of the ship indicates the speed of the ship for a given propulsor(s) input in calm water, i.e., without waves.
    
    In the case of level-1 vulnerability criterion, a ship is considered to be nonvulnerable to the broaching failure mode if
    \begin{equation}
        L \ge 200~[\mathrm{m}]
        \label{eq:L1_L}
    \end{equation}
    or
    \begin{equation}
        \mathrm{Fn}\le 0.3
        \label{eq:L1_Fn}
    \end{equation}
    Here, $L~[\mathrm{m}]$ is the length of the ship, as defined in paragraph 2.12 of the introduction part of the 2008 IS Code, and
    $\mathrm{Fn}$ is the Froude number, which is a speed-length ratio based on inertial and gravitational forces defined as
    \begin{equation}
        \mathrm{Fn} = \frac{V_{\rm{s}}}{\sqrt{Lg}}
    \end{equation}
    Here, $V_{\rm{s}}$ represents the service speed of the ship and $g$ is the gravitational acceleration.
    Eq.~\ref{eq:L1_L} corresponds to the fact that a ship can only surf-ride a steep wave of a length comparable to the length of the ship and that a long, steep wave rarely occurs in the ocean.
    Eq.~\ref{eq:L1_Fn} corresponds to the lowest surf-riding threshold in the case of conventional ships \cite{Imo2022-uf}.

    In the level-2 vulnerability criterion, vulnerability is judged based on the occurrence probability of surf-riding.
    In ocean engineering, sea states are generally considered from two perspectives: short- and long-term.
    In the level-2 vulnerability criterion, the short-term sea-state statistics are characterized by the significant wave height and zero-crossing wave period, and
    the long-term sea-state statistics are given as a joint frequency table of the significant wave height and the zero-crossing wave period (wave scatter table). 
    Once the spectral density of the wave elevation is achieved for the short-term sea state,
    the probability density function of the local regular wave is determined according to Longuet--Higgins' theory \cite{Longuet-Higgins1983-xa}.
    Then, the critical nominal speed of the ship is calculated for each local regular wave, and
    the occurrence probability of surf-riding in the sea state is estimated as 
    the occurrence probability of the local regular waves, where the critical nominal speed of the ship is less than the service speed of the ship.
    Therefore, in the case of the level-2 vulnerability criterion, a ship is considered to be nonvulnerable to the broaching failure mode if
    \begin{equation}
        C \le R_\mathrm{SR}
    \end{equation}
    where
    \begin{gather}
        C = \sum_{H_\mathrm{S}} \sum_{T_\mathrm{Z}} \left[ \mathit{W2} \left( H_\mathrm{S} , T_\mathrm{Z} \right) \sum_{i=0}^{N_\lambda} \sum_{j=0}^{N_a} w_{ij}\left( H_\mathrm{S} , T_\mathrm{Z} \right) \mathit{C2}_{ij} \right] \\
        R_\mathrm{SR} = 0.005
    \end{gather}
    Here, 
    $C$: a value corresponding to the occurrence probability of surf-riding,
    $\mathit{W2} \left( H_\mathrm{S} , T_\mathrm{Z} \right)$: probability of short-term sea states based on the wave scatter table (Table~\ref{tab:WaveScatter}) \cite{Iacs2001-bg},
    $H_\mathrm{S}$: significant wave height [m] that is specified in Table~\ref{tab:WaveScatter},
    $T_\mathrm{Z}$: zero-crossing wave period [s] that is specified in Table~\ref{tab:WaveScatter},
    $w_{ij}$: joint probability density function of a local wave with the wave steepness and wave-length to ship-length ratio \cite{Longuet-Higgins1983-xa},
    $N_a$: number of discretization for the wave-length to ship-length ratio of the local regular wave,
    $N_\lambda$: number of discretization for the steepness of the local regular wave, and
    $\mathit{C2}_{ij}$: the coefficient of the occurrence of surf-riding.
    The joint probability-density function of a local wave, $w_{ij}$, is calculated as
    \begin{equation}
        \begin{aligned}
            &w_{i j\left( H_\mathrm{S} , T_\mathrm{Z} \right)}=4 \frac{\sqrt{g}}{\pi v} \frac{L^{\frac{5}{2}} T_{01}}{H_s{}^3} s_j{}^2 r_i{}^{\frac{3}{2}} \frac{\sqrt{1+v^2}}{1+\sqrt{1+v^2}}~\Delta r~\Delta s~ \cdot\\
            &\exp \left\{-2\left(\frac{L r_i s_j}{H_s}\right)^2\left[1+\frac{1}{v^2}\left(1-\sqrt{\frac{g T_{01}{}^2}{2 \pi r_i L}}\right)^2\right]\right\}
        \end{aligned}
        \label{eq:wij}
    \end{equation}
    Here, $\nu$:represents the band parameter, $0.425$;
    $T_{01}$: is the mean wave period, $1.086 T_\mathrm{Z}$, 
    $s_j$:is the wave steepness of the local regular wave varying from 0.03 to 0.15 with the increment of $\Delta s = 0.0012$, and
    $r_i$ is the wave-length to ship-length ratio of the local regular wave varying from 1.0 to 3.0 with the increment of $\Delta r = 0.025$.
    The coefficient on the occurrence of surf-riding $\mathit{C2}_{ij}$ is defined as
    \begin{equation}
        \label{eq:cij}
        \mathit{C2}_{i j} \quad=\left\{\begin{array}{lll}
        1 & \text { if } & \mathrm{Fn} >    \mathrm{Fn}_{\mathrm{cr}}\left(r_j, s_i\right) \\
        0 & \text { if } & \mathrm{Fn} \leq \mathrm{Fn}_{\mathrm{cr}}\left(r_j, s_i\right)
        \end{array}\right.
    \end{equation}
    \begin{equation}
        \mathrm{Fn}_{\mathrm{cr}}= \frac{u_{\mathrm{cr}}}{\sqrt{L g}}
    \end{equation}
    Here, 
    $\mathrm{Fn}_{\mathrm{cr}}$ represents the critical Froude number
    in regular waves with respect to steepness $s_j$ and wave-length to ship-length ratio $r_i$ and
    $u_{\mathrm{cr}}$ is the critical nominal speed of the ship [m/s].
    The critical nominal speed of the ship, $u_{\mathrm{cr}}$, is estimated based on Melnikov' method (Section ~\ref{section_PhasePortrait}).
    If the ship uses propeller(s) as the propulsor(s),
    the relationship between the critical nominal speed of the ship and the critical number of revolutions of the propeller(s) $n_\mathrm{cr}$ is defined as
    \begin{equation}
        T_e\left( u_\mathrm{cr}; n_\mathrm{cr} \right) - R\left( u_\mathrm{cr} \right) = 0
        \label{eq:criticalspeedandrev}
    \end{equation}
    where 
    $T_e\left( u; n \right)$ is the propeller thrust in calm waters [N],
    $R\left( u \right)$ is the ship resistance in calm waters [N],
    $u$ is the speed of the ship [m/s], and 
    $n$ is the number of revolutions of the propeller(s) [1/s].
    The propeller thrust and ship resistance should be balanced in calm waters.
    In addition, the propeller thrust in calm waters, $T_e\left( u; n \right)$, is approximated as
    \begin{equation}
        \begin{aligned}
            &T_e\left( u; n \right)\\
            &=\left(1-t_\mathrm{p}\right) \rho n^2 D_\mathrm{p}^4\left\{\kappa_0+\kappa_1 J\left( u, n \right)+\kappa_2 \left[J\left( u, n \right)\right]^2\right\} \\
            &=\tau_0 n^2 + \tau_1 u n + \tau_2 u^2
        \end{aligned}
    \end{equation}
    where
    \begin{align}
        J\left( u, n \right)&= \frac{ u \left( 1 - w_\mathrm{p} \right)}{n D_\mathrm{p}} \\
        \tau_0&=\kappa_0 \left( 1-t_\mathrm{p} \right) \rho D_\mathrm{p}{}^4 \\
        \tau_1&=\kappa_1 \left( 1-t_\mathrm{p} \right) \left( 1-w_\mathrm{p} \right) \rho D_\mathrm{p}{}^3 \\
        \tau_2&=\kappa_2 \left( 1-t_\mathrm{p} \right) \left( 1-w_\mathrm{p} \right)^2 \rho D_\mathrm{p}{}^2
    \end{align}
    where 
    $J$ is the advance ratio;
    $t_\mathrm{p}$ is the approximate thrust deduction factor;
    $w_\mathrm{p}$ is the approximate wake fraction;
    $\kappa_0$, $\kappa_1$, and $\kappa_2$ are the approximation coefficients for the approximated propeller thrust coefficient in calm water; and
    $\tau_0$, $\tau_1$, and $\tau_2$ are the approximation coefficients for the approximated propeller thrust coefficient in calm waters as functions of $u$ and $n$
    The ship resistance in calm waters, $R\left( u \right)$, is approximated by the quintic polynomial as
    \begin{equation}
        R(u) = r_1 u+r_2 u^2+r_3 u^3+r_4 u^4+r_5 u^5,
    \end{equation}
    where $r_1$, $r_2$, $r_3$, $r_4$, and $r_5$ represent the regression coefficients for the ship resistance in calm waters.
    In addition, the amplitude of the wave surging force by the local regular wave $f_{ij}$ [N] is estimated as
    \begin{equation}
        f_{i j}=\rho g k_i \frac{H_{i j}}{2} \sqrt{F c_i{}^2+F s_i{}^2}
    \end{equation}
    where
    \begin{equation}
        H_{i j}= s_j r_i L
    \end{equation}
    \begin{equation}
        k_i=\frac{2 \pi}{r_i L}
    \end{equation}
    \begin{equation}
        Fc_i=\sum_{m=1}^N \delta x_m~S\left(x_m\right) \sin \left(k_i~x_m\right) \exp \left[-\frac{1}{2} k_i~d\left(x_m\right)\right]
    \end{equation}
    \begin{equation}
        Fs_i=\sum_{m=1}^N \delta x_m~S\left(x_m\right) \cos \left(k_i~x_m\right) \exp \left[-\frac{1}{2} k_i~d\left(x_m\right)\right]
    \end{equation}
    where $H_{i j}$ is the height of the local regular wave [m],
    $k_i$ is the wave number of the local regular wave [1/m],
    $Fc_i$ and $Fs_i$ represent the parts of the Froude--Krylov component of the wave surging force [m],
    $m$ is the index of a station, $N$ is the number of stations,
    and $x_m$ is the longitudinal distance from the midship to the station $m$ [m] (positive to the bow section). Furthermore,
    $\delta x_m$ is the length of the ship strip associated with the station $m$ [m],
    $S\left(x_m\right)$ represents the area of the submerged portion of the ship at station $m$ in calm waters [$\mathrm{m}^2$], and 
    $d\left(x_m\right)$: the draft at the station $m$ in calm water [m]. Therefore, the critical number of revolutions of the propeller(s) $n_\mathrm{cr}\left( s_j,r_i \right)$ is directly derived from Melnikov's method shown in Section~\ref{sec:Melnikov_Hamiltonian_part} as
    \begin{equation}
        \begin{aligned}
            &2 \pi \frac{T_e\left(c_i; n_{c r}\right)-R\left(c_i\right)}{f_{i j}}+8 a_0 n_{c r}+8 a_1-4 \pi a_2+\frac{64}{3} a_3\\
            &-12 \pi a_4+\frac{1024}{15} a_5=0
        \end{aligned}
        \label{eq:MelnikovIMO}
    \end{equation}
    where
    \begin{align}
        a_0&=-\frac{\tau_1}{\sqrt{f_{i j}~k_i \left(M+M_x\right)}} \\
        a_1&=\frac{r_1+2 r_2 c_i+3 r_3 c_i^2+4 r_4 c_i^3+5 r_5 c_i^4-2 \tau_2 c_i}{\sqrt{f_{i j}~k_i \left(M+M_x\right)}}\\
        a_2&=\frac{r_2+3 r_3 c_i+6 r_4 c_i^2+10 r_5 c_i^3-\tau_2}{k_i \left(M+M_x\right)}\\
        a_3&=\frac{r_3+4 r_4 c_i+10 r_5 c_i^2}{\sqrt{k_i^3\left(M+M_x\right)^3}} \sqrt{f_{i j}}\\
        a_4&=\frac{r_4+5 r_5 c_i}{k_i^2\left(M+M_x\right)^2} f_{i j}\\
        a_5&=\frac{r_5}{\sqrt{k_i^5\left(M+M_x\right)^5}} \sqrt{f_{i j}{}^3}\\
        c_i&= \sqrt{\frac{g}{k_i}}
    \end{align}
    Here,
    $M$: the mass of the ship [kg],
    $M_x$: the added mass of the ship in surge [kg], and
    $c_i$: the wave celerity of the local regular wave [m/s].
    Sakai et al. has shown that Eq.~\ref{eq:MelnikovIMO} is a quadratic equation of $n_\mathrm{cr}$ and 
    that the larger solution is appropriate for the critical number of revolutions of the propeller(s) if Eq.~\ref{eq:MelnikovIMO} has two real solutions \cite{Sakai2017-ns}.
    The critical nominal speed of the ship $u_{\mathrm{cr}}$ is obtained by Eq.~\ref{eq:criticalspeedandrev}.

    Notably, Eq.~\ref{eq:MelnikovIMO} is completely identical with the estimation formula of the surf-riding threshold (Eq.~\ref{eq:Melnikov_5th} in this paper), as shown by Maki~\cite{maki2010surfriding}.
    \begin{table*}[h]
        \caption{Wave scatter table used in the SGISC \cite{Iacs2001-bg}. The total number of occurrences is $10^6$.}
        \centering
        \begin{tabular}{@{}S|@{~}S@{~}S@{~}S@{~}S@{~}S@{~}S@{~}S@{~}S@{~}S@{~}S@{~}S@{~}S@{~}S@{~}S@{~}S@{~}S@{}}
            \hline \hline
            \multicolumn{1}{c|@{~}}{\backslashbox{$H_\mathrm{S}~[\mathrm{m}]$}{$T_\mathrm{Z}~[\mathrm{s}]$}} & 
                  3.5 & 4.5 & 5.5 & 6.5 & 7.5 & 8.5 & 9.5 & 10.5 \\ \hline
            0.5 & 1.3 & 133.7 & 865.6 & 1186.0 & 634.2 & 186.3 & 36.9 & 5.6 \\ 
            1.5 & 0.0 & 29.3 & 986.0 & 4976.0 & 7738.0 & 5569.7 & 2375.7 & 703.5 \\ 
            2.5 & 0.0 & 2.2 & 197.5 & 2158.8 & 6230.0 & 7449.5 & 4860.4 & 2066.0 \\ 
            3.5 & 0.0 & 0.2 & 34.9 & 695.5 & 3226.5 & 5675.0 & 5099.1 & 2838.0 \\ 
            4.5 & 0.0 & 0.0 & 6.0 & 196.1 & 1354.3 & 3288.5 & 3857.5 & 2685.5 \\ 
            5.5 & 0.0 & 0.0 & 1.0 & 51.0 & 498.4 & 1602.9 & 2372.7 & 2008.3 \\ 
            6.5 & 0.0 & 0.0 & 0.2 & 12.6 & 167.0 & 690.3 & 1257.9 & 1268.6 \\ 
            7.5 & 0.0 & 0.0 & 0.0 & 3.0 & 52.1 & 270.1 & 594.4 & 703.2 \\ 
            8.5 & 0.0 & 0.0 & 0.0 & 0.7 & 15.4 & 97.9 & 255.9 & 350.6 \\ 
            9.5 & 0.0 & 0.0 & 0.0 & 0.2 & 4.3 & 33.2 & 101.9 & 159.9 \\ 
            10.5 & 0.0 & 0.0 & 0.0 & 0.0 & 1.2 & 10.7 & 37.9 & 67.5 \\ 
            11.5 & 0.0 & 0.0 & 0.0 & 0.0 & 0.3 & 3.3 & 13.3 & 26.6 \\ 
            12.5 & 0.0 & 0.0 & 0.0 & 0.0 & 0.1 & 1.0 & 4.4 & 9.9 \\ 
            13.5 & 0.0 & 0.0 & 0.0 & 0.0 & 0.0 & 0.3 & 1.4 & 3.5 \\ 
            14.5 & 0.0 & 0.0 & 0.0 & 0.0 & 0.0 & 0.1 & 0.4 & 1.2 \\ 
            15.5 & 0.0 & 0.0 & 0.0 & 0.0 & 0.0 & 0.0 & 0.1 & 0.4 \\ 
            16.5 & 0.0 & 0.0 & 0.0 & 0.0 & 0.0 & 0.0 & 0.0 & 0.1 \\ \hline 
            \multicolumn{1}{c}{~} & ~ & ~ & ~ & ~ & ~ & ~ & ~ & ~ \\ \hline \hline
            \multicolumn{1}{c|@{~}}{\backslashbox{$H_\mathrm{S}~[\mathrm{m}]$}{$T_\mathrm{Z}~[\mathrm{s}]$}} &
                 11.5 & 12.5 & 13.5 & 14.5 & 15.5 & 16.5 & 17.5 & 18.5 \\ \hline
            0.5 & 0.7 & 0.1 & 0.0 & 0.0 & 0.0 & 0.0 & 0.0 & 0.0 \\ 
            1.5 & 160.7 & 30.5 & 5.1 & 0.8 & 0.1 & 0.0 & 0.0 & 0.0 \\ 
            2.5 & 644.5 & 160.2 & 33.7 & 6.3 & 1.1 & 0.2 & 0.0 & 0.0 \\ 
            3.5 & 1114.1 & 337.7 & 84.3 & 18.2 & 3.5 & 0.6 & 0.1 & 0.0 \\ 
            4.5 & 1275.2 & 455.1 & 130.9 & 31.9 & 6.9 & 1.3 & 0.2 & 0.0 \\ 
            5.5 & 1126.0 & 463.6 & 150.9 & 41.0 & 9.7 & 2.1 & 0.4 & 0.1 \\ 
            6.5 & 825.9 & 386.8 & 140.8 & 42.2 & 10.9 & 2.5 & 0.5 & 0.1 \\ 
            7.5 & 524.9 & 276.7 & 111.7 & 36.7 & 10.2 & 2.5 & 0.6 & 0.1 \\ 
            8.5 & 296.9 & 174.6 & 77.6 & 27.7 & 8.4 & 2.2 & 0.5 & 0.1 \\ 
            9.5 & 152.2 & 99.2 & 48.3 & 18.7 & 6.1 & 1.7 & 0.4 & 0.1 \\ 
            10.5 & 71.7 & 51.5 & 27.3 & 11.4 & 4.0 & 1.2 & 0.3 & 0.1 \\ 
            11.5 & 31.4 & 24.7 & 14.2 & 6.4 & 2.4 & 0.7 & 0.2 & 0.1 \\ 
            12.5 & 12.8 & 11.0 & 6.8 & 3.3 & 1.3 & 0.4 & 0.1 & 0.0 \\ 
            13.5 & 5.0 & 4.6 & 3.1 & 1.6 & 0.7 & 0.2 & 0.1 & 0.0 \\ 
            14.5 & 1.8 & 1.8 & 1.3 & 0.7 & 0.3 & 0.1 & 0.0 & 0.0 \\ 
            15.5 & 0.6 & 0.7 & 0.5 & 0.3 & 0.1 & 0.1 & 0.0 & 0.0 \\ 
            16.5 & 0.2 & 0.2 & 0.2 & 0.1 & 0.1 & 0.0 & 0.0 & 0.0 \\  \hline
        \end{tabular}
        \label{tab:WaveScatter}
    \end{table*}
    
\section{Conclusion}
This paper presented a comprehensive review of all analytical formulae proposed so far for estimating the surf-riding threshold. The equation of motion governing the surge motion of vessels in following seas is not directly solvable because of its nonlinearity. To address this difficulty, approximate solution methods were used along with a numerical solution approach. The results show that the proposed methods were able to capture the qualitative trend of the surfing threshold. Among them, the approximate estimation method based on Melnikov's method could quantitatively capture the trend of the surf-riding threshold. Finally, the interconnection between the prediction formula rooted in Melnikov's method and its relevance with IMO's SGISC was comprehensively explained.

\begin{acknowledgements}
This work was supported by a Grant-in-Aid for Scientific Research from the Japan Society for Promotion of Science (JSPS KAKENHI Grant Numbers JP22H01701, JP21K14362). The authors are also thankful to Enago (www.enago .jp) for reviewing the English language.
\end{acknowledgements}

% could have direct or potential influence or impact bias on 
\section*{Conflict of interest}

    The authors declare that they have no conflict of interest.

% \bibliographystyle{spmpsci}      % mathematics and physical sciences
\bibliographystyle{spphys}       % APS-like style for physics
\bibliography{paper.bib}   % name your BibTeX database

\begin{thebibliography}{10}
\providecommand{\url}[1]{{#1}}
\providecommand{\urlprefix}{URL }
\expandafter\ifx\csname urlstyle\endcsname\relax
  \providecommand{\doi}[1]{DOI \discretionary{}{}{}#1}\else
  \providecommand{\doi}{DOI \discretionary{}{}{}\begingroup
  \urlstyle{rm}\Url}\fi

\bibitem{saunders1965}
H.E. Saunders, The preliminary hydrodynamics design of a motorboat, Chapter 77
  of Hydrodynamics in Ship Design  (1965)

\bibitem{nicholson1974}
K.~Nicholson, Some parametric model experiment to investigate broaching-to.
\newblock in \emph{Proc Int Symp Dynamics of Marine Vehicle and Structure,
  1974} (1974)

\bibitem{kan1990capsizing}
M.~Kan, T.~Saruta, H.~Taguchi, M.~Yasuno, Capsizing of a ship in quatering seas
  part 1. model experiments on mechanism of capsizing, Journal of the Society
  of Naval Architects of Japan \textbf{1990}(167), 81 (1990)

\bibitem{spyrou1996dynamic}
K.~Spyrou, Dynamic instability in quartering seas—part ii: analysis of ship
  roll and capsize for broaching, Journal of Ship Research \textbf{40}(04), 326
  (1996)

\bibitem{spyrou1997dynamic}
K.~Spyrou, Dynamic instability in quartering seas—part iii: nonlinear effects
  on periodic motions, Journal of Ship Research \textbf{41}(03), 210 (1997)

\bibitem{maki2009bifurcation}
A.~Maki, N.~Umeda, Bifurcation and chaos in yaw motion of a ship at lower speed
  in waves and its prevention using optimal control.
\newblock in \emph{Proceedings of the 10th International Conference on
  stability of ships and ocean vehicles} (2009), pp. 429--440

\bibitem{motora1982mechanism}
S.~Motora, M.~Fujino, T.~Fuwa, On the mechanism of broaching-to phenomena.
\newblock in \emph{Proceedings of the 2nd international conference on stability
  of ships and ocean vehicles.} (1982), pp. 535--550

\bibitem{Grim1951}
O.~Grim, \emph{Das Schiff in von achtern auflaufender See} (Springer Berlin
  Heidelberg, Berlin, Heidelberg, 1951), pp. 264--287

\bibitem{Makov1969}
Y.~Makov, Some results of theoretical analysis of surfriding in following seas,
  T krylov Soc \textbf{126} (1969)

\bibitem{Ananiev1966}
D.~Ananiev, On surf-riding in following seas, T Krylov Soc \textbf{13}, 169
  (1966)

\bibitem{Kan_1990_surfriding}
M.~Kan, A guideline to avoid the dangerous surf-riding, Proceedings of the 4th
  international conference on stability of ships and ocean vehicles pp. 90--97
  (1990)

\bibitem{Umeda1990SR_irregular}
N.~Umeda, Probabilistic study on surf-riding of a ship in irregular following
  seas, Proceedings of the 4th International Conference on Stability of Ships
  and Ocean Vehicles pp. 336--343 (1990)

\bibitem{Imo1995}
{IMO}.
\newblock {MSC.-Circ}. 707 - guidance to the master for avoiding dangerous
  situations in following and quatering seas (1995)

\bibitem{holmes1980averaging}
P.J. Holmes, Averaging and chaotic motions in forced oscillations, SIAM Journal
  on Applied Mathematics \textbf{38}(1), 65 (1980)

\bibitem{spyrou2006}
K.~Spyrou, Asymmetric surging of ships in following seas and its repercussions
  for safety, Nonlinear Dynamics \textbf{43}, 149 (2006)

\bibitem{maki2010melnikov}
A.~Maki, N.~Umeda, T.~Ueta, Melnikov integral formula for beam sea roll motion
  utilizing a non-hamiltonian exact heteroclinic orbit, Journal of marine
  science and technology \textbf{15}(1), 102 (2010)

\bibitem{spyrou2001exact}
K.~Spyrou, Exact analytical solutions for asymmetric surging and surf-riding.
\newblock in \emph{Proceeding of 5th International Workshop on Stability and
  Operational Safety of Ships, University of Trieste, Trieste}, vol.~4 (2001),
  vol.~4, p.~3

\bibitem{maki2010surfriding}
A.~Maki, N.~Umeda, M.~Renilson, T.~Ueta, Analytical formulae for predicting the
  surf-riding threshold for a ship in following seas, Journal of marine science
  and technology \textbf{15}, 218 (2010)

\bibitem{maki2014melnikov}
A.~Maki, N.~Umeda, T.~Ueta, Melnikov integral formula for beam sea roll motion
  utilizing a non-hamiltonian exact heteroclinic orbit: analytic extension and
  numerical validation, Journal of Marine Science and Technology
  \textbf{19}(3), 257 (2014)

\bibitem{Umeda1983}
N.~Umeda, On the surf-riding of a ship, Journal of the Society of Naval
  Architects of Japan \textbf{1982}(152), 192 (1983)

\bibitem{umeda1990SR_regular}
N.~Umeda, T.~Kohyama, Surf-riding of a ship in regular seas (in japanese),
  Journal of the Kansai Society of Naval Architects, Japan \textbf{213}, 63
  (1990)

\bibitem{kuznetsov1998elements}
Y.A. Kuznetsov, I.A. Kuznetsov, Y.~Kuznetsov, \emph{Elements of applied
  bifurcation theory}, vol. 112 (Springer, 1998)

\bibitem{kawakami1981separatrix}
H.~Kawakami, Calculating method of separatrix loop of saddle-typed fixed
  points, J Inst Electron Inf Commun Eng, A \textbf{64}(10), 860 (1981)

\bibitem{doedel1990numerical}
E.J. Doedel, M.J. Friedman, Numerical computation of heteroclinic orbits,
  Continuation Techniques and Bifurcation Problems pp. 155--170 (1990)

\bibitem{friedman1991numerical}
M.J. Friedman, E.J. Doedel, Numerical computation and continuation of invariant
  manifolds connecting fixed points, SIAM Journal on Numerical Analysis
  \textbf{28}(3), 789 (1991)

\bibitem{maki2014surfriding}
A.~Maki, N.~Umeda, M.~Renilson, T.~Ueta, Analytical methods to predict the
  surf-riding threshold and the wave-blocking threshold in astern seas, Journal
  of Marine Science and Technology \textbf{19}, 415 (2014)

\bibitem{nagumo1962active}
J.~Nagumo, S.~Arimoto, S.~Yoshizawa, An active pulse transmission line
  simulating nerve axon, Proceedings of the IRE \textbf{50}(10), 2061 (1962)

\bibitem{mckean1970nagumo}
H.~McKean~Jr, Nagumo's equation, Advances in mathematics \textbf{4}(3), 209
  (1970)

\bibitem{komuro1992bifurcation}
M.~Komuro, Bifurcation equations of continuous piecewise-linear vector fields,
  Japan journal of industrial and applied mathematics \textbf{9}, 269 (1992)

\bibitem{endo1993piecewise}
T.~Endo, L.O. Chua, Piecewise-linear analysis of high-damping chaotic
  phase-locked loops using melnikov's method, IEEE Transactions on Circuits and
  Systems I: Fundamental Theory and Applications \textbf{40}(11), 801 (1993)

\bibitem{Belenky1993}
V.L. Belenky, A capsizing probability computation method, Journal of Ship
  Research \textbf{37}(3), 200 (1993)

\bibitem{Belenky1994}
V.L. Belenky, Piecewise linear method for the probabilistic stability
  assessment for ship in a seaway, 5th International Conference on Stability of
  Ship and Ocean Vehicle \textbf{5}(2), 1 (1994)

\bibitem{Guckenheimer_Holmes1983}
J.~Guckenheimer, P.~Holmes, \emph{Nonlinear oscillations, dynamical systems,
  and bifurcations of vector fields}, vol.~42 (Springer Science \& Business
  Media, 2013)

\bibitem{Greenspan_Holmes1984}
B.~Greenspan, P.~Holmes, Repeated resonance and homoclinic bifurcation in a
  periodically forced family of oscillators, SIAM journal on mathematical
  analysis \textbf{15}(1), 69 (1984)

\bibitem{Imo2020-bp}
{IMO}.
\newblock {MSC.1-Circ}. 1627 - interim guidelines on the second generation
  intact stability criteria (2020)

\bibitem{Imo2022-uf}
{IMO}.
\newblock {MSC.1-Circ}. 1652 - explanatory notes to the interim guidelines on
  second generation intact stability criteria (2022)

\bibitem{maki2016surfriding}
A.~Maki, Y.~Miyauchi, Prediction methods for the surf-riding threshold and the
  wave-blocking threshold based on melnikov’s method, Journal of Marine
  Science and Technology \textbf{21}(2), 179 (2016)

\bibitem{salam1987mel}
F.M. Salam, The mel’nikov technique for highly dissipative systems, SIAM
  Journal on Applied Mathematics \textbf{47}(2), 232 (1987)

\bibitem{wu_mccue2008application}
W.~Wu, L.~McCue, Application of the extended melnikov's method for
  single-degree-of-freedom vessel roll motion, Ocean Engineering
  \textbf{35}(17-18), 1739 (2008)

\bibitem{maki2022NOLTA}
A.~Maki, Y.~Miino, N.~Umeda, M.~Sakai, T.~Ueta, H.~Kawakami, Nonlinear dynamics
  of ship capsizing at sea, Nonlinear Theory and Its Applications, IEICE
  \textbf{13}(1), 2 (2022)

\bibitem{Longuet-Higgins1983-xa}
M.S. Longuet-Higgins, On the joint distribution of wave periods and amplitudes
  in a random wave field, Proc. R. Soc. Lond. A Math. Phys. Sci.
  \textbf{389}(1797), 241 (1983)

\bibitem{Iacs2001-bg}
{IACS}.
\newblock Recommendation {N}o. 34 ({C}orr. 1 {N}ov. 2001): Standard wave data
  (2001)

\bibitem{Sakai2017-ns}
M.~Sakai, A.~Maki, T.~Murakami, N.~Umeda, Analytical solution of critical speed
  for {Surf-Riding} in the light of melnikov analysis.
\newblock in \emph{Proceedings of the Japan Society of Naval Architects and
  Ocean Engineers}, vol.~24 (The Japan Society of Naval Architects and Ocean
  Engineers, 2017), vol.~24, pp. 311--314

\end{thebibliography}

\end{document}